\documentclass[12pt,a4paper,twoside,final,notitlepage, leqno]{article}
\usepackage[english]{babel}
\usepackage[T1]{fontenc}
\usepackage{epsfig, graphicx, amssymb}
\usepackage{amsmath,amsthm,epsfig,amsfonts}
\usepackage{float}
\usepackage{color}

\newcommand{\monitem}{ \smallskip \noindent $\bullet$ \quad  }
\newcommand{\moneq}{\vspace*{-7pt} \begin{equation} \displaystyle }
\newcommand{\moneqstar}{\vspace*{-6pt} \begin{equation*} \displaystyle }
\newcommand{\monendstar}{\vspace*{-6pt} \end{equation*}   }
\newcommand{\monend}{\vspace*{-7pt} \end{equation}   }
\newcommand{\moneqarraystar}{ \begin{eqnarray*} \displaystyle }
\newcommand{\monendarraystar}{ \end{eqnarray*}   }

\newcommand{\R}{\mathbb{R}}

\newcommand{\N}{\mathbb{N}}

\definecolor{vertfonce}{rgb}{0.0, 0.5, 0.0}

\newcommand{\rotV}{\mathbf{rot}}

\def \proof {\noindent \textbf{Proof}}
\def \monitem {\noindent{$\bullet~$}}
\newtheorem{theorem}{Theorem}
\newtheorem{lemme}{Lemma}
\newtheorem{defi}{Definition}
\newtheorem{prop}{Proposition}
\newtheorem{rem}{Remark}

\usepackage{fancyhdr}
\renewcommand{\headrulewidth}{0pt}

\begin{document}

\fancypagestyle{plain}{ \fancyfoot{} \renewcommand{\footrulewidth}{0pt}}
\fancypagestyle{plain}{ \fancyhead{} \renewcommand{\headrulewidth}{0pt}}

~

  \vskip 2.1 cm

\centerline {\bf \LARGE Discrete harmonics}

 \bigskip 

\centerline {\bf \LARGE for stream function-vorticity Stokes problem} 

\bigskip  \bigskip \bigskip

\centerline { \large   Fran\c{c}ois Dubois, Michel Sala\"un, St\'ephanie Salmon}

\smallskip  \bigskip

\centerline { \it  \small
  IAT / CNAM UPRES EA n$^{\mathrm{o}} 2140$}

\centerline { \it  \small   15, rue Marat, F-78210 Saint-Cyr-L'Ecole.}


\bigskip  \bigskip

\centerline {April 17, 2001 
{\footnote {\rm  \small $\,$  Published in {\it Numerische Mathematik},  volume 92, number 4,
  pages 711-742, October 2002. Edition 12  December 2024.}}}

 \bigskip \bigskip

\bigskip  \bigskip
\noindent {\bf \large Abstract}

\noindent
We consider the bidimensional Stokes problem for incompressible fluids in stream function-vorticity. For this problem, 
the classical finite elements method of degree one converges only in $\mathcal{O}(\sqrt{h})$ for the $L^2-$norm 
of the vorticity. We propose to use harmonic functions to approach the vorticity along the boundary. Discrete harmonics 
are functions that are used in practice to derive a new numerical method. We prove that we obtain with this 
numerical scheme an error of order $\mathcal{O}(h)$ for the $L^2-$norm of the vorticity.


 \bigskip \bigskip
{\bf{AMS (MOS) classification}} : $65$N$30$.

\noindent

\newpage



\section{Introduction}

\fancyhead[EC]{\sc{Fran\c{c}ois Dubois, Michel Sala\"un, St\'ephanie Salmon}}
\fancyhead[OC]{\sc{Discrete harmonics for stream function-vorticity Stokes problem}}
\fancyfoot[C]{\oldstylenums{\thepage}}

\subsection{Motivations}

Let $\Omega$ be a bounded convex polygonal open set in $\R^2$. We denote $\Gamma = \partial \Omega$ 
the boundary of $\Omega$.
Modelization of the equilibrium of an incompressible and viscous fluid leads to the Navier-Stokes problem 
(Landau-Lifschitz \cite{LL}). 
But, when viscosity is sufficiently important or velocity of the fluid sufficiently small, 
we can neglect the convection terms and we obtain the stationary Stokes problem which is,  
in primitive variables (velocity $\mathbf{u}$ and pressure $p$)~:
\begin{eqnarray*}
\left\{
\begin{array}{ccccc}
-\nu \Delta \mathbf{u} + \nabla p & = & \mathbf{f} & \mbox{ in } & \Omega \\
\mbox{div}~\mathbf{u} & = & 0 & \mbox{ in } & \Omega \\
\mathbf{u} & = & 0 & \mbox{ on } & \Gamma
\end{array}
\right.
\end{eqnarray*}
where $\nu$ is the kinematic viscosity and $\mathbf{f}$ a field of given external forces. \\

As velocity is divergence free, this two dimensional problem is often rewritten with 
stream function and vorticity. Velocity is the curl of some stream function 
and vorticity is the curl of the velocity. Let us remark that, classically, 
in stream function-vorticity formulation, the velocity is assumed to be zero 
on the whole boundary~: for more general boundary conditions, we refer to \cite{BCMP}, 
\cite{G88} and \cite{BL}. Then, a usual way of discretizing the Stokes problem, according to this 
formulation, is to choose a finite element method with polynomials of degree one. The stream function and the 
vorticity are searched in $H^1(\Omega)$, but it is well-known that the problem is not mathematically 
well-posed (see for example Girault-Raviart \cite{GR}). Nevertheless, it can be shown that 
the scheme is convergent (Scholz \cite{Sch}, \cite{GR}), the convergence for the $L^2-$norm 
of the vorticity is in $\sqrt{h}$, where $h$ is the maximum diameter of triangles in the 
triangulation and $h^{1-\varepsilon}$ for the $H^1-$norm of the stream function. 
There exists a variational formulation of the Stokes problem which is well-posed but difficult to use 
numerically (Bernardi-Girault-Maday \cite{BGM}). We propose in this paper to study this 
well-posed formulation in order to obtain a numerical scheme convergent in at least $\mathcal{O}(h)$ 
for the $L^2-$norm of the vorticity. In the last part of the paper, we also give results of 
numerical experiments.

\subsection{Notation}

We shall consider the following spaces (see for example \cite{Adams}). 
We denote $\mathcal{D}(\Omega)$ the space of all indefinitely differentiable functions from $\Omega$ to 
$\R$ with compact support and $L^2(\Omega)$ the space of all classes of functions which are square integrable. 
For any integer $m \geq 0$ and any real number $p$ satisfying $1 \leq p \leq \infty$,
$W^{m,p}(\Omega)$ consists of those functions $v \in L^p(\Omega)$ for which all partial derivatives $\partial^{\alpha}v$ 
(in the distribution sense) with $|\alpha| \leq m$ belong to the space $L^p(\Omega)$ :\\
$$W^{1,2}(\Omega) = H^1(\Omega) = \left\{\varphi \in L^2(\Omega)~,~\forall i \in \left\{1,2\right\}~,~
\frac{\partial \varphi}{\partial x_i} \in L^2(\Omega) \right\}$$
$$W^{2,2}(\Omega) = H^2(\Omega) = \left\{ \varphi \in H^1(\Omega)~,~\forall i \in \left\{1,2\right\}~,~
\frac{\partial \varphi}{\partial x_i} \in H^1(\Omega) \right\}~~~.$$
Symbols $\parallel {\scriptstyle{\bullet}} \parallel_{_{j,\Omega}}$ 
(respectively $|{\scriptstyle{\bullet}}|_{_{j,\Omega}}$) 
denote norms (respectively semi-norms) in Sobolev spaces $H^j(\Omega)$, $j=1,2.$ 
The Sobolev space $H^1_0(\Omega)$ (respectively $H^2_0(\Omega)$) is the closure of $\mathcal{D}(\Omega)$
in the sense of the norm $\parallel {\scriptstyle{\bullet}} \parallel_{_{1,\Omega}}$ 
(respectively $\parallel {\scriptstyle{\bullet}} \parallel_{_{2,\Omega}}$).
In the sequel, $({\scriptstyle{\bullet}},{\scriptstyle{\bullet}})$ denotes the standard inner product in $L^2(\Omega)$ 
and $\langle {\scriptstyle{\bullet}},{\scriptstyle{\bullet}}\rangle_{_{-1,1}}$, the duality product between $H^1_0(\Omega)$ 
and its topological dual space $H^{-1}(\Omega)$.
Finally, $\gamma$ shall denote the trace operator from $H^1(\Omega)$ onto $H^{1/2}(\Gamma)$, or 
from $H^2(\Omega)$ onto $H^{3/2}(\Gamma)$ (see Lions-Magenes \cite{LM}).

\section{The Stokes problem}

Let $\mathbf{f}$ be a field of given forces and suppose that $\mathbf{f}$ belongs to $(L^2(\Omega))^2$ 
and $\displaystyle{\mbox{rot}~\mathbf{f} = \frac{\partial f_1}{\partial x_2} - \frac{\partial f_2}{\partial x_1}}$ 
belongs to $L^2(\Omega).$
The stationary Stokes problem consists in finding a stream function $\psi$ and a vorticity
field $\omega$ solutions of~:
\begin{eqnarray}
\omega + \Delta \psi                                     & = & 0 ~~~~~~ \mbox{ in } \Omega     \label{E1} \\
-\Delta \omega                                           & = & \mbox{rot}~\mathbf{f} ~ \mbox{ in } \Omega     \label{E2} \\
\psi                                                     & = & 0 ~~~~~~ \mbox{ on } \Gamma     \label{E3} \\
\displaystyle{\frac{\partial \psi}{\partial \mathbf{n}}} & = & 0 ~~~~~~ \mbox{ on } \Gamma~~~. \label{E4} 
\end{eqnarray}
Indeed, it consists in finding a velocity field $\mathbf{u}$ that is divergence free and can
be written by means of a stream function $\psi$~:
$$\mathbf{u} = \rotV~\psi~~~\mbox{in}~\Omega~~~.$$
Equation (\ref{E1}) means that the vorticity is the curl of the velocity. Equation (\ref{E2})
is the equilibrium equation for a viscous fluid, with kinematic viscosity equal to $1$, and
where convection terms are neglected. The boundary conditions (\ref{E3}) and (\ref{E4}) are
consequences of~;
$$\mathbf{u} = 0~~~\mbox{on}~\Gamma~~~.$$
For more details about this problem, we refer to \cite{GR}.\\

It is natural to discretize the problem (\ref{E1})-(\ref{E4}) with a linear and continuous
finite element method.
As $\Omega$ is supposed to be polygonal, we exactly cover it with a mesh $\mathcal{T}$ composed of 
triangular elements.
The mesh $\mathcal{T}$ is supposed to be regular (in the sense defined in Ciarlet \cite{Ci}). The set  
$H^1_{_{\mathcal{T}}}$ denotes the space of continuous functions defined on
$\overline{\Omega}$, polynomial of degree $1$ in each triangle of $\mathcal{T}$ and
$H^1_{_{0,\mathcal{T}}} = H^1_{_{\mathcal{T}}} \cap H^1_0(\Omega)$~:
$$H^1_{_{\mathcal{T}}} = \left\{ \varphi \in \mathcal{C}^0(\overline{\Omega})~,~\forall K \in
\mathcal{E}_{_{\mathcal{T}}}~,~\varphi_{|_K} \in {\mathbb{P}}^1(K) \right\}~~~,$$                           
where $\mathcal{E}_{_{\mathcal{T}}}$ is the set of triangles in $\mathcal{T}$ and $ {\mathbb{P}}^1$ the    
space of all polynomials of total degree $\leq 1.$
The discretization of the problem (\ref{E1})-(\ref{E4}) consists in finding
$\omega_{_{\mathcal{T}}}$ and $\psi_{_{\mathcal{T}}}$~:
\begin{eqnarray}
& \psi_{_{\mathcal{T}}} \in H^1_{_{0,\mathcal{T}}} (\mbox{boundary condition (\ref{E3}) is
then satisfied}) & \label{E5}\\
& \omega_{_{\mathcal{T}}} \in H^1_{_{\mathcal{T}}} & \label{E6}
\end{eqnarray}
in the following way. We first multiply (\ref{E1}) by a scalar function 
$\varphi \in H^1_{_{\mathcal{T}}}$ and integrate by parts $(\Delta \psi,\varphi)$. 
Taking the boundary condition (\ref{E4}) into account, we obtain~:
\begin{eqnarray}
(\omega_{_{\mathcal{T}}},\varphi) - (\nabla \psi_{_{\mathcal{T}}},\nabla \varphi) = 0~~,~~
\forall \varphi \in H^1_{_{\mathcal{T}}}~~~. \label{E7}
\end{eqnarray}
Then, we multiply (\ref{E2}) by a test function $\xi$ in $H^1_{_{0,\mathcal{T}}}$, 
and integrate by parts both sides. As $\xi$ vanishes on the boundary, we obtain~:
\begin{eqnarray}
(\nabla \omega_{_{\mathcal{T}}}, \nabla \xi) = (\mathbf{f},\rotV~\xi)~~,~~\forall \xi \in
H^1_{_{0,\mathcal{T}}}~~~. \label{E8}
\end{eqnarray}
This formulation (\ref{E5})-(\ref{E8}) has been studied for a long time (see for example
Ciarlet-Raviart \cite{CR}, Glowinski-Pironneau \cite{GP}, \cite{GR} among others).
But results are not optimal.
On the one hand, the continuous formulation associated is not well-posed for any
$\mathbf{f}$ in $(L^2(\Omega))^2.$ Indeed, the Stokes problem could be seen
as a biharmonic problem for the stream function~: 
$$\psi \in H^2_0(\Omega)$$
$$\Delta^2 \psi = \mbox{rot}~\mathbf{f}~~~.$$
In a variational form, this problem can be rewritten as follows~:
$$(\Delta \psi,\Delta \mu) = (\mathbf{f},\rotV~\mu)~~,~~\forall \mu \in H^2_0(\Omega)~~~.$$
This problem is well-posed in $H^2_0(\Omega)$ (see Ciarlet \cite{Ci}) but
vorticity cannot be more regular than square integrable. On the other hand, it is impossible 
to obtain polynomial functions of degree $1$ which really satisfy both
boundary conditions (\ref{E3}) and (\ref{E4}) without being identically zero on elements $K$
in contact with $\Gamma$ (\emph{i.e} $\partial K \cap \Gamma$ nor empty, nor reduced to one point).
Numerically, the Neumann condition (\ref{E4}) is never satisfied because weakly written.
Moreover, error estimates derived for this scheme are not optimal. The bound is in $h_{_{\mathcal{T}}}^{1/2}$
for the $L^2$-norm of vorticity, where $h_{_{\mathcal{T}}}$ is the maximum diameter of the elements
of the triangulation \cite{Sch}, \cite{GR} (see Figure {\ref{convns} at the end of the article). 
Nevertheless, it is possible to stabilize the numerical scheme by adding jumps at interfaces 
and prove convergence, see \cite{AB}.

A different weak formulation was proposed by Bernardi, Girault and Maday
\cite{BGM}, who introduced the space~: 
$$M(\Omega) = \left\{ \varphi \in L^2(\Omega)~,~\Delta \varphi \in H^{-1}(\Omega) \right\}~~~.$$ 
It consists in finding $(\omega,\psi)$ in $M(\Omega) \times H^1_0(\Omega)$ with the following method. 
We test the first equation (\ref{E1}) with a function $\varphi \in M(\Omega)$ and the second one 
(\ref{E2}) with a function $\xi \in H^1_0(\Omega)$. It is then necessary to integrate 
twice by parts the term $\langle \Delta \psi,\varphi\rangle_{_{-1,1}}$ in order to include 
boundary conditions (\ref{E3}) and (\ref{E4}). Then, we get~: 
\begin{equation} \label{P}
\left\{
\begin{array}{llll}
(\omega,\varphi)~+~\langle \Delta \varphi,\psi\rangle_{_{-1,1}} & = & 0                             & \forall \varphi \in M(\Omega)\\
\langle -\Delta \omega,\xi\rangle_{_{-1,1}}                     & = & (\mathbf{f},\mathbf{rot}~\xi) & \forall \xi \in H^1_0(\Omega)~~~.
\end{array}
\right.
\end{equation}
This formulation (\ref{P}) leads to a well-posed problem \cite{BGM}.

\begin{rem} \label{rem1}
The space $M(\Omega)$ is a Hilbert space for the norm:
$$\parallel \varphi \parallel_{_M} = \sqrt{\parallel \varphi \parallel^2_{_{0,\Omega}}
+ \parallel \Delta \varphi \parallel^2_{_{-1,\Omega}}} \,.$$
\end{rem}

\begin{prop}{} \label{rem2} 
The Sobolev space $H^1(\Omega)$ is contained in $M(\Omega)$ with continuous imbedding. 
Moreover, we have~:
$$\forall \varphi \in H^1(\Omega)~~,~~\parallel \varphi \parallel_{_M}~\leq~\parallel \varphi \parallel_{_{1,\Omega}}~~~.$$ 
Finally, if $\varphi$ belongs to $M(\Omega) \cap H^1_0(\Omega) = H^1_0(\Omega)$,
the M-norm is equivalent to the $H^1$-norm.
\end{prop}

\proof \\
First, we remark that the laplacian of a function in $H^1(\Omega)$ belongs to $H^{-1}(\Omega)$.
Then, to prove the continuity of the imbedding, we write the definition of the $H^{-1}$-norm. For all
$\varphi \in H^1(\Omega)$, we have~:
\begin{equation} \label{norme.H-1}
\parallel \Delta \varphi \parallel_{_{-1,\Omega}}
~=~\sup_{\chi \in H^1_0(\Omega)} \frac{\langle \Delta \varphi,\chi\rangle_{_{-1,1}}}{|\chi |_{_{1,\Omega}}}
~=~\sup_{\chi \in H^1_0(\Omega)} \frac{-(\nabla \varphi,\nabla \chi)}{\parallel \nabla \chi \parallel_{_{0,\Omega}}}~~~,
\end{equation}
because $\varphi$ is in $H^1(\Omega)$ and $\chi_{|\Gamma} = 0$. Using the Cauchy-Schwarz inequality, we obtain~:
$$\parallel \Delta \varphi \parallel_{_{-1,\Omega}}~\leq~\parallel \nabla \varphi \parallel_{_{0,\Omega}}~~~,$$
and then~:
$$\parallel \varphi \parallel_{_M}~\leq~\parallel \varphi \parallel_{_{1,\Omega}}~~~.$$
Finally, if $\varphi$ belongs to $H^1_0(\Omega)$, it suffices to take 
$\chi = \varphi$ in (\ref{norme.H-1}) and we obtain~:  
$\parallel \varphi \parallel_{_M}~=~\parallel \varphi \parallel_{_{1,\Omega}}$. So the M-norm 
is equivalent to the $H^1_0$-norm. $\hfill \blacksquare$\\

The analysis of (\ref{P}), conducted in \cite{BGM}, is straightforward with arguments of Brezzi \cite{Br}. We
introduce the kernel $\mathcal{H}(\Omega)$ of the bilinear form $\langle \Delta .,.\rangle_{_{-1,1}}$~:
$$\mathcal{H}(\Omega) = \left\{ \varphi \in M(\Omega)~,~\langle \Delta \varphi,\xi\rangle_{_{-1,1}} = 0~~~
\forall \xi \in H^1_0(\Omega) \right\}~~~.$$

\begin{prop}{Characterization of the space $\mathcal{H}(\Omega)$.} \label{prop1} \\
Let $\Omega$ be a bounded Lipschitz open subset of $\R^2$. Then, we have~:
$$\mathcal{H}(\Omega) = \left\{ \varphi \in L^2(\Omega)~,~\Delta \varphi = 0 \mathrm{~in~} L^2(\Omega) \right\}~~~,$$
and there exists a trace operator, still denoted by $\gamma$, from $\mathcal{H}(\Omega)$ on $H^{-1/2}(\Gamma)$. 
\end{prop}

\proof \\
\monitem It is obvious that the space $\mathcal{H}(\Omega)$ is a subset of the space~:
$$D(\Delta,L^2(\Omega))~=~\left\{ \varphi \in L^2(\Omega)~,~\Delta \varphi \in L^2(\Omega) \right\}~~~.$$
This space was extensively studied by Grisvard \cite{Gri} when $\Omega$ is a bounded polygonal open subset of $\R^2$, 
who recalled that the space $H^2(\Omega)$ is dense in $D(\Delta,L^2(\Omega))$. For completeness of our study, 
the proof of this result is given below. \\
\monitem It is based on the following property~:
a subspace $\mathcal{S}$ of a Hilbert space $M$ is dense in $M$ if and only if every element
of $M'$ that vanishes on $\mathcal{S}$ also vanishes on $M$. 
Let $\widehat{u}$ belong to $(D(\Delta,L^2(\Omega)))'$. As it is a Hilbert space, the Riesz theorem 
proves that there exists a function, denoted by $u$, in $D(\Delta,L^2(\Omega))$ such that~:
$$\langle \widehat{u} , \varphi \rangle_{_{(D(\Delta,L^2(\Omega)))',D(\Delta,L^2(\Omega))}}
~=~(u , \varphi)_{_{D(\Delta,L^2(\Omega))}}~~,~~\forall~\varphi \in D(\Delta,L^2(\Omega))~~~.$$
Using the expression of the $D(\Delta,L^2(\Omega))$ scalar product, we have for all $\varphi$ in 
$D(\Delta,L^2(\Omega))$~:
$$(u , \varphi)_{_{D(\Delta,L^2(\Omega))}}~=~(u , \varphi )~+~(\Delta u , \Delta \varphi )~~~.$$
We suppose now that $\widehat{u}$ vanishes on $H^2(\Omega)$. Then, we obtain~:
$$(u , \varphi )~+~(\Delta u , \Delta \varphi )~=~0~,~\forall \varphi \in H^2(\Omega)~~~.$$
Let us now introduce $\widetilde{u}$ and $\psi$ the extensions to $\R^2$, 
by zero outside $\Omega$, of functions $u$ and $\Delta u$. Moreover, let us 
notice that, for all function $\widetilde{\varphi} \in \mathcal{D}(\R^2)$, its restriction on $\Omega$, 
say $\varphi$, belongs to $H^2(\Omega)$. Then, the above formula leads to the following relations~:
\begin{eqnarray*}
0 & = & (u , \varphi )~+~(\Delta u , \Delta \varphi )\\
  & = & (\widetilde{u} , \varphi )~+~(\psi , \Delta \varphi )\\
  & = & (\widetilde{u},\widetilde{\varphi})_{_{L^2(\scriptsize{\R^2})}} + 
        (\psi,\Delta \widetilde{\varphi})_{_{L^2(\scriptsize{\R}^2)}}
	~~~\forall \widetilde{\varphi} \in \mathcal{D}(\R^2)~~~.
\end{eqnarray*}
This equality implies that in the distributions sense~:
$$\widetilde{u} + \Delta \psi = 0 ~\mathrm{ in } ~\mathcal{D'}(\R^2)~~~.$$
Thus, $\Delta~\psi$ belongs to $L^2(\R^2)$ since $\widetilde{u}$ belongs to $L^2(\R^2)$. Moreover, 
$\psi$ belongs to $L^2(\R^2)$ as, by definition, $\Delta u$ belongs to $L^2(\Omega)$. Then, we deduce 
that $\psi$ belongs to $H^2(\R^2)$ (use the characterization of the space $H^2(\R^2)$ with the Fourier 
transform). As $\psi$ is identically zero outside $\Omega$, we deduce that $\psi$ belongs to
$H^2_0(\Omega)$ (see \cite{GR}). Using the density of $\mathcal{D}(\Omega)$ in $H^2_0(\Omega)$, 
we introduce a sequence $(\psi_k)_{k \geq 1}$ of $\mathcal{D}(\Omega)$ that tends to $\psi$ in 
$H^2_0(\Omega)$. Then, we have the following relations~:
\begin{eqnarray*}
\psi_k        & \stackrel{k \longrightarrow \infty}{\longrightarrow} & \psi \mathrm{~in~} L^2(\Omega)~~~,\\
\Delta \psi_k & \stackrel{k \longrightarrow \infty}{\longrightarrow} & \Delta \psi = - u \mathrm{~in~} L^2(\Omega)~~~.
\end{eqnarray*}
Using the above convergences, we obtain for all $\varphi \in D(\Delta,L^2(\Omega))$~:
\begin{eqnarray*}
(u , \varphi)_{_{D(\Delta,L^2(\Omega))}} & = & (u , \varphi )~+~(\Delta u , \Delta \varphi ) \\
& = & (u , \varphi )~+~(\psi , \Delta \varphi ) \\
& = & \lim_{k \longrightarrow \infty} [~- (\Delta \psi_k , \varphi)~+~(\psi_k , \Delta \varphi)~] \\
& = & \lim_{k \longrightarrow \infty} [~- \langle \varphi , \Delta \psi_k \rangle_{_{\mathcal{D'}(\Omega),\mathcal{D}(\Omega)}}
      ~+~\langle \Delta \varphi , \psi_k \rangle_{_{\mathcal{D'}(\Omega),\mathcal{D}(\Omega)}}~]~=~0~~~.
\end{eqnarray*}
Then we obtain~: $(u , \varphi)_{_{D(\Delta,L^2(\Omega))}}~=~0$ for any $\varphi$ in $D(\Delta,L^2(\Omega))$
which implies that $u$ is zero and finishes to prove the density of $H^2(\Omega)$ in $D(\Delta,L^2(\Omega))$.

\monitem Now, let us recall the Green's formula, which is valid in any bounded Lipschitz domain 
(see Ne\u{c}as \cite{Nec})~:
$$\int_{\Omega} (u~\Delta v~-~v~\Delta u) \, \mbox{d}x~=~\int_{\Gamma} (\gamma u~\frac{\partial v}{\partial n}
            ~-~\gamma v~\frac{\partial u}{\partial n}) \, \mbox{d}\Gamma~~~,$$
for all $u$ and $v$ in $H^2(\Omega)$; $\displaystyle{\frac{\partial u}{\partial n}}$ is the normal 
derivative of $u$ along $\Gamma$. Then, due to the density of $H^2(\Omega)$ in $D(\Delta,L^2(\Omega))$, 
the Green's formula permits to define a trace operator from $\mathcal{H}(\Omega)$ on $H^{-1/2}(\Gamma)$. 
$\hfill \blacksquare$\\

So, when we restrict the first equation of (\ref{P}) to functions in $\mathcal{H}(\Omega)$, this new
problem is well-posed according to Lax-Milgram's lemma \cite{LMi}, as~:

\begin{prop}\label{prop2}
The $L^2-$scalar product~: \\
$$\mathcal{H}(\Omega) \times \mathcal{H}(\Omega) \ni (\omega,\varphi) \longmapsto 
\displaystyle{\int_{\Omega} \omega~\varphi~\mathrm{d}x} \in \R$$
is elliptic on $\mathcal{H}(\Omega)$.
\end{prop}

\proof \\
For $\omega$ in $\mathcal{H}(\Omega)$, it is obvious that~:
$\parallel \omega \parallel^2_{_M}~=~\parallel \omega  \parallel^2_{_{0,\Omega}}$~. $\hfill \blacksquare$

\begin{prop} {Decomposition of space $M(\Omega)$.} \label{prop0} \\
We have~: 
$$ M(\Omega) = H^1_0(\Omega) \oplus \mathcal{H}(\Omega)~~~.$$
\end{prop}

\proof \\
We split $\varphi \in M(\Omega)$ into two parts~: $\varphi = \varphi^0 + \varphi^{\Delta}$. 
One the one hand, since $\Delta \varphi \in H^{-1}(\Omega)$, the component $\varphi^0$ is 
uniquely defined in $H^1_0(\Omega)$ by the Dirichlet problem~:
$$
\left\{
\begin{array}{lll}
\Delta \varphi^0 & = & \Delta \varphi~\mbox{ in } \Omega \\
\gamma \varphi^0 & = & 0~~~~ \mbox{ on } \Gamma~~~.
\end{array}
\right.
$$
On the other hand, we define the function $\varphi^{_{\Delta}}$ by $\varphi^{_{\Delta}} = \varphi - \varphi^0$.  
Then it satisfies~: 
$$\Delta \varphi^{\Delta} = 0~~~~\mbox{ in } \Omega~~~.$$ 
So $\varphi^{_{\Delta}}$ belongs to $\mathcal{H}(\Omega)$ and moreover~:
$$\gamma \varphi^{_{\Delta}} = \gamma \varphi~~~\mbox{ on } \Gamma~~~,$$
with $\gamma \varphi$ well defined in $H^{-1/2}(\Gamma)$ (see Proposition \ref{prop1}).
$\hfill \blacksquare$

\section{Harmonic discretization of Stokes problem}

\subsection{Harmonic lifting}

We recall that $\Omega$ being polygonal allows to cover it entirely with a mesh
$\mathcal{T}$. We denote $\mathcal{E}_{_{\mathcal{T}}}$ the set of triangles in $\mathcal{T}$, 
$\mathcal{S}_{_{\mathcal{T}}}$ the set of vertices of $\mathcal{T}$ and 
$h_{_{\mathcal{T}}} = \sup \left( \mbox{diam}~K~;~K\in \mathcal{E}_{_{\mathcal{T}}} \right).$ 

\begin{defi}{Family $\mathcal{U}_{\sigma}$ of regular meshes.} \label{usigma}\\
We suppose that $\mathcal{T}$ belongs to the set 
$\mathcal{U}_{\sigma}$ of triangulations satisfying :
$$\exists~\sigma > 0~,~\forall K \in \mathcal{E}_{_{\mathcal{T}}}~~,~~\frac{h_K}{\rho_K} \leq \sigma~~~,$$  
where $h_K = \mathrm{diam}~K$ and $\rho_K$ is the diameter of the circle inscribed in $K$.  
\end{defi}

\noindent We define the space $\Gamma_{\mathcal{T}}$ of all traces of functions in
$H^1_{_{\mathcal{T}}}$. Dimension of $\Gamma_{\mathcal{T}}$ is exactly the number 
$NS(\mathcal{T},\Gamma)$ of vertices of the mesh $\mathcal{T}$ on the boundary $\Gamma$.
\begin{eqnarray*}
\Gamma_{\mathcal{T}} = \left\{ 
\begin{array}{r}
\lambda : \partial \Omega = \Gamma \longrightarrow \R~,~\lambda \mathrm{~continuous~on~} \Gamma \\
\lambda \mathrm{~linear~on~each~edge~of~the~mesh~} 
\end{array}
\right\}~~~.
\end{eqnarray*}

\begin{defi}{Harmonic lifting.}
$$\mathcal{H}_{_{\mathcal{T},\infty}} = \left\{ \varphi \in H^1(\Omega)~,~\exists \lambda \in
\Gamma_{\mathcal{T}}~,~\gamma \varphi = \lambda~,~\Delta \varphi = 0 \mathrm{~in~}
H^{-1}(\Omega) \right\}~~~.$$
\end{defi}

\noindent The discrete space $\mathcal{H}_{_{\mathcal{T},\infty}}$ is finite-dimensional and it is
clear that $\mathcal{H}_{_{\mathcal{T},\infty}}$ and $\Gamma_{\mathcal{T}}$ have the same dimension.   
We define the following harmonic lifting operator $Z_{\infty}$~: 
$$Z_{\infty} : \Gamma_{_{\mathcal{T}}} \ni \lambda \longmapsto 
Z_{\infty}(\lambda) \in \mathcal{H}_{_{\mathcal{T},\infty}}$$
by the Dirichlet problem~: 
\begin{eqnarray} \label{zinfini}
\left\{
\begin{array}{rcl}
\Delta Z_{\infty}(\lambda) & = & 0~~\mbox{ in } \Omega\\
\gamma Z_{\infty}(\lambda) & = & \lambda ~~\mbox{ on } \Gamma~~~.
\end{array}
\right.
\end{eqnarray}

\begin{rem}\label{rem3}
It is useful to remark that $\mathcal{H}_{_{\mathcal{T},\infty}}$ is included in
$\mathcal{H}(\Omega) \subset M(\Omega).$
\end{rem}

\begin{defi}{Harmonic discretization of $H^1(\Omega)$.}\\
We set~:
\begin{equation} \label{split}
H^{1,\infty}_{_{\mathcal{T}}} = H^1_{_{0,\mathcal{T}}} \oplus \mathcal{H}_{_{\mathcal{T},\infty}}~~~.
\end{equation}
\end{defi}

\noindent The dimension of $H^{1,\infty}_{_{\mathcal{T}}}$ is the same as $H^1_{_{\mathcal{T}}}$ but, on
the boundary, we replace piecewise linear continuous functions by harmonic functions, whose 
support is spread over the entire domain but whose values decrease rapidly.

\begin{lemme}{Property of $\mathcal{H}_{_{\mathcal{T},\infty}}$.} \label{prop3b}
\begin{equation} \label{gradgrad}
\forall \varphi \in \mathcal{H}_{_{\mathcal{T},\infty}}~,~\forall \xi \in H^1_{_{0,\mathcal{T}}}~,~
(\nabla \varphi, \nabla \xi) = 0~~~.
\end{equation}
\end{lemme}

\begin{rem}
The splitting $\varphi = \varphi^0 + \varphi^{\Delta}$ described in Proposition \ref{prop0} 
remains valid for functions in $H^{1,\infty}_{_{\mathcal{T}}}$. Indeed, for all function 
$\varphi$ of $H^{1,\infty}_{_{\mathcal{T}}}$, $\gamma \varphi$ belongs to $\Gamma_{_{\mathcal{T}}}$.  
Then $\varphi^{\Delta} = Z_{\infty}(\gamma\varphi)$ is in $\mathcal{H}_{_{\mathcal{T},\infty}}$, and,  
by construction, $\varphi^0 = \varphi - \varphi^{\Delta}$ belongs to space $H^1_{_{0,\mathcal{T}}}$.
\end{rem}

\subsection{Formulation}

We propose the following harmonic variational formulation of the Stokes problem~:
\begin{eqnarray} \label{E9-11}
& & \mbox{Find } \psi_{_{\mathcal{T}}} \in H^1_{_{0,\mathcal{T}}}~\mbox{and}~\omega_{_{\mathcal{T}}} \in
H^{1,\infty}_{_{\mathcal{T}}} \mbox{ such that : } \nonumber \\
& & \left\{
\begin{array}{lcll}
(\omega_{_{\mathcal{T}}},\varphi) - (\nabla \psi_{_{\mathcal{T}}},\nabla \varphi) & = & 0 
&\forall~\varphi \in H^{1,\infty}_{_{\mathcal{T}}} \\
(\nabla \omega_{_{\mathcal{T}}}, \nabla \xi)                                      & = & (\mathbf{f},\rotV~\xi) 
& \forall~\xi \in H^1_{_{0,\mathcal{T}}}~~~.
\end{array}
\right.
\end{eqnarray}
As $H^{1,\infty}_{_{\mathcal{T}}}$ is included in $H^1(\Omega)$, $\nabla \omega_{_{\mathcal{T}}}$ is well
defined. The method, proposed here, is a natural discretization of problem (\ref{P}) since, 
as it was said in Remark \ref{rem3}, $H^{1,\infty}_{_{\mathcal{T}}}$ is a subset of $M(\Omega)$.

\begin{prop}{Existence and uniqueness of a solution to problem (\ref{E9-11}).} \label{prop3}\\
If $\mathbf{f} \in (L^2(\Omega))^2$, the problem (\ref{E9-11}) has a unique solution
$(\psi_{_{\mathcal{T}}},\omega_{_{\mathcal{T}}})$ in the space $H^1_{_{0,\mathcal{T}}} \times
H^{1,\infty}_{_{\mathcal{T}}}$ which depends continuously on the datum $\mathbf{f}$. There exists 
a strictly positive constant $C$ independent of the mesh such that~:
\begin{equation} \label{ex1}
\parallel \omega_{_{\mathcal{T}}} \parallel_{_M} + \parallel \nabla \psi_{_{\mathcal{T}}} \parallel_{_{0,\Omega}}
\leq C \parallel \mathbf{f} \parallel_{_{0,\Omega}}.
\end{equation}
\end{prop}

\proof \\
\monitem We split all functions in $H^{1,\infty}_{_{\mathcal{T}}}$ into two parts as explicited in (\ref{split}) 
and we rewrite the problem (\ref{E9-11}) using such a splitting for $\omega_{_{\mathcal{T}}} = 
\omega^0_{_{\mathcal{T}}} + \omega^{\Delta}_{_{\mathcal{T}}}$ and for
the test function $\varphi = \varphi^0 + \varphi^{\Delta}$. Discrete functions $\omega^0_{_{\mathcal{T}}}$ 
and $\varphi^0$ belong to the space $H^1_{_{0,\mathcal{T}}}$ while  
$\omega^{\Delta}_{_{\mathcal{T}}}$ and $\varphi^{\Delta}$ lie in the space $\mathcal{H}_{_{\mathcal{T},\infty}}$~:
\begin{eqnarray*}
\left\{
\begin{array}{llll}
\psi_{_{\mathcal{T}}} \in H^1_{_{0,\mathcal{T}}}~,~ 
\omega^0_{_{\mathcal{T}}} \in H^1_{_{0,\mathcal{T}}}~,~ 
\omega^{\Delta}_{_{\mathcal{T}}} \in \mathcal{H}_{_{\mathcal{T},\infty}} \\
(\omega^0_{_{\mathcal{T}}},\varphi^0) + (\omega^{\Delta}_{_{\mathcal{T}}},\varphi^0)
- (\nabla \psi_{_{\mathcal{T}}},\nabla \varphi^0) & = & 0 & \forall \varphi^0 \in
H^1_{_{0,\mathcal{T}}}\\
(\omega^0_{_{\mathcal{T}}},\varphi^{\Delta}) +
(\omega^{\Delta}_{_{\mathcal{T}}},\varphi^{\Delta})
- (\nabla \psi_{_{\mathcal{T}}},\nabla \varphi^{\Delta}) & = & 0 & 
\forall \varphi^{\Delta} \in \mathcal{H}_{_{\mathcal{T},\infty}} \\
(\nabla \omega^0_{_{\mathcal{T}}}, \nabla \xi) + (\nabla \omega^{\Delta}_{_{\mathcal{T}}},
\nabla \xi) & = & (\mathbf{f},\rotV~\xi) & \forall \xi \in H^1_{_{0,\mathcal{T}}}~~~.
\end{array}
\right.
\end{eqnarray*}
Due to Lemma \ref{prop3b}, we obtain~:
\begin{eqnarray}
(\omega^0_{_{\mathcal{T}}},\varphi^0) + (\omega^{\Delta}_{_{\mathcal{T}}},\varphi^0)
- (\nabla \psi_{_{\mathcal{T}}},\nabla \varphi^0) & = & 0~~~~~~~~~~~~\forall \varphi^0 \in
H^1_{_{0,\mathcal{T}}} \label{E12}\\
(\omega^0_{_{\mathcal{T}}},\varphi^{\Delta}) +
(\omega^{\Delta}_{_{\mathcal{T}}},\varphi^{\Delta})
& = & 0~~~~~~~~~~~~\forall \varphi^{\Delta} \in \mathcal{H}_{_{\mathcal{T},\infty}} \label{E13}\\
(\nabla \omega^0_{_{\mathcal{T}}}, \nabla \xi) & = & (\mathbf{f},\rotV~\xi)~~~\forall \xi \in
H^1_{_{0,\mathcal{T}}}.
\label{E14}
\end{eqnarray}
The way of studying problem (\ref{E12})-(\ref{E14}) follows ideas of \cite{GP}, Achdou, Glowinski, 
Pironneau \cite{AGP} and Ruas \cite{Ru1}.\\
\monitem Problem (\ref{E14}) is well-posed according to Lax-Milgram's lemma, and Poincar\'e's
inequality. As $\Omega$ is bounded, if $C_p$ denotes the Poincar\'e constant~:
$$\forall \chi \in H^1_0(\Omega)~~,~~\parallel \chi \parallel_{_{0,\Omega}}~\leq~C_p~| \chi |_{_{1,\Omega}}~~~,$$ 
there exists a unique $\omega^0_{_{\mathcal{T}}} \in H^1_{_{0,\mathcal{T}}}$ satisfying (take
$\xi = \omega^0_{_{\mathcal{T}}}$ in (\ref{E14}))~:
$$\parallel \nabla \omega^0_{_{\mathcal{T}}} \parallel^2_{_{0,\Omega}}~=~
\parallel \rotV~\omega^0_{_{\mathcal{T}}} \parallel^2_{_{0,\Omega}}~ 
\leq~\parallel \mathbf{f} \parallel_{_{0,\Omega}}~\parallel \rotV~\omega^0_{_{\mathcal{T}}}\parallel_{_{0,\Omega}}~~~.$$
The $M$-norm of $\omega^0_{_{\mathcal{T}}}$ is equivalent to its $H^1$-norm, as 
$\omega^0_{_{\mathcal{T}}} \in H^1_0(\Omega)$ (see Proposition \ref{rem2}), and we obtain~:
\begin{equation} \label{ex11}
\parallel \omega^0_{_{\mathcal{T}}} \parallel_{_M}~\leq~\sqrt{1+C^2_p} \parallel \mathbf{f} \parallel_{_{0,\Omega}}~~~. 
\end{equation}
\monitem We study now equation (\ref{E13}) where the function $\omega^0_{_{\mathcal{T}}}$ is given.  
Using the ellipticity of the $L^2$-scalar product on $\mathcal{H}_{_{\mathcal{T},\infty}}$ 
(see Proposition \ref{prop2}), the problem~: 
$$(\omega^{\Delta}_{_{\mathcal{T}}},\varphi^{\Delta})~=~- (\omega^0_{_{\mathcal{T}}},\varphi^{\Delta})
~~,~~\forall \varphi^{\Delta} \in \mathcal{H}_{_{\mathcal{T},\infty}}~~~,$$
has a unique solution $\omega^{\Delta}_{_{\mathcal{T}}}$, which verifies 
(take $\varphi^{\Delta} = \omega^{\Delta}_{_{\mathcal{T}}}$)~:
\begin{equation} \label{ex12}
\parallel \omega^{\Delta}_{_{\mathcal{T}}} \parallel_{_M}~=~
\parallel \omega^{\Delta}_{_{\mathcal{T}}} \parallel_{_{0,\Omega}}~\leq~
\parallel \omega^0_{_{\mathcal{T}}} \parallel_{_{0,\Omega}}~\leq~
\sqrt{1+C^2_p} \parallel \mathbf{f} \parallel_{_{0,\Omega}}~~~.
\end{equation}
\monitem Finally, equation (\ref{E12}) is formally identical to (\ref{E14}), so there exists
(Lax-Milgram's lemma and Poincar\'e's inequality) $\psi_{_{\mathcal{T}}} \in H^1_{_{0,\mathcal{T}}}$ such that~:
$$(\nabla \psi_{_{\mathcal{T}}},\nabla \varphi^0)~=~(\omega^0_{_{\mathcal{T}}},\varphi^0)
~+~(\omega^{\Delta}_{_{\mathcal{T}}},\varphi^0)~~,~~\forall \varphi^0 \in H^1_{_{0,\mathcal{T}}}~~~.$$
Then, using (\ref{ex11}) and (\ref{ex12}), we deduce that~:
\begin{equation} \label{ex13}
\parallel \nabla \psi_{_{\mathcal{T}}} \parallel_{_{0,\Omega}}~\leq~
C_p \parallel \omega^0_{_{\mathcal{T}}} \parallel_{_{0,\Omega}}~+~ 
C_p \parallel \omega^{\Delta}_{_{\mathcal{T}}} \parallel_{_{0,\Omega}}~\leq~
C'\parallel \mathbf{f} \parallel_{_{0,\Omega}}~~~.
\end{equation}
Combining (\ref{ex11}),(\ref{ex12}) and (\ref{ex13}), we obtain (\ref{ex1}).$\hfill \blacksquare$

\begin{prop}{Regularity of harmonic liftings.} \label{prop8}\\
Let $\varepsilon$ be a strictly positive real number. We have~:
\begin{equation} \label{reg}
\mathcal{H}_{_{\mathcal{T},\infty}} \subset H^{2-\varepsilon}(\Omega)~~,~~\forall \varepsilon > 0~~~. 
\end{equation}
\end{prop}
\begin{figure}    [H]  \centering
\centerline  {\includegraphics[width=.50\textwidth]   {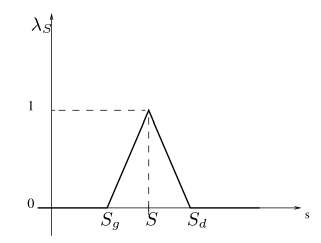}}
\caption{A basis function.}
\label{lambda}
\end{figure}

\proof \\
For a given vertex $S$ of $\mathcal{T}$ on the boundary $\Gamma$ ($S \in \mathcal{S}_{_{\mathcal{T}}} 
\cap \Gamma$), a basis function $\lambda_S$ 
of $\Gamma_{_{\mathcal{T}}}$ is represented on Figure \ref{lambda} where $s$ is the curvilinear 
abscissa along $\Gamma$. 
Thus, the gradient of such a function is a combination of Heaviside functions, so belongs to
$H^{1/2-\varepsilon}(\Gamma)$ for any $\varepsilon > 0$ but not to $H^{1/2}(\Gamma)$ \cite{Adams}, 
\cite{LM}. A basis function of space $\mathcal{H}_{_{\mathcal{T},\infty}}$ is defined by finding 
$\varphi_S : \Omega \longrightarrow \R$ such that $\varphi_S \in H^1(\Omega),~
\Delta \varphi_S = 0, \gamma \varphi_S = \lambda_S$ \emph{i.e} $\varphi_S = Z_{\infty}(\lambda_S)$. 
Then, the imbedding (\ref{reg}) is thus a consequence of regularity of the laplacian operator on a convex domain
\cite{Gri}. $\hfill \blacksquare$\\

To obtain error estimates and thus convergence, we need a stability
result (Proposition \ref{prop4}) and an interpolation error (Proposition \ref{prop5}).

\subsection{Interpolation error}

Let $\Pi_{_{\mathcal{T}}} : H^2(\Omega) \longrightarrow H^1_{_{\mathcal{T}}}$ be the Lagrange 
interpolation operator associated with the mesh $\mathcal{T}$, and 
$\Pi^{\Gamma}_{_{\mathcal{T}}} : H^1(\Gamma) \longrightarrow \Gamma_{_{\mathcal{T}}}$ 
the analogous interpolation operator on the boundary $\Gamma$. Let $h_{_{\mathcal{T}}}$ be the maximum diameter of 
triangles in $\mathcal{T}.$

\begin{prop}{Property of the interpolation operator $\Pi^{\Gamma}_{_{\mathcal{T}}}$ on the boundary.} 
\label{proppit}\\
If $\mathcal{T}$ belongs to a family $\mathcal{U}_{\sigma}$ of regular meshes, 
the interpolation operator $\Pi^{\Gamma}_{_{\mathcal{T}}}$ is defined from $H^{3/2}(\Gamma)$ 
on $H^{1/2}(\Gamma) \cap \Gamma_{_{\mathcal{T}}}$ and we have~: 
\begin{equation} \label{pit} 
\parallel \mu - \Pi^{\Gamma}_{_{\mathcal{T}}}\mu \parallel_{_{1/2,\Gamma}} 
~\leq~C(\sigma)~h_{_{\mathcal{T}}}~\parallel \mu \parallel_{_{3/2,\Gamma}}~~,~~\forall \mu \in H^{3/2}(\Gamma)~~~.
\end{equation}
\end{prop}

\noindent To prove this result, we need the following theorem \cite{LM}~:

\begin{theorem}{Interpolation between Sobolev spaces.} \label{Iss}\\
Let $s_i$ and $t_i$ be two couples of positive real numbers for $i = 0$ or $i=1$ and 
$p \in \R$ such that $1 \leq p < \infty$. Let $\Pi$ be an operator of
$\mathcal{L}(W^{s_0,p}(\Omega);W^{t_0,p}(\Omega)) \cap
\mathcal{L}(W^{s_1,p}(\Omega);W^{t_1,p}(\Omega))$ ($\mathcal{L}$ 
is the space of all linear and continuous functions). Then,
for all $\theta \in \R \cap [0;1]$, operator $\Pi$ belongs to the interpolate space \\
$\mathcal{L}_{\theta} = \mathcal{L}(W^{(1-\theta)s_0+\theta
s_1,p}(\Omega);W^{(1-\theta)t_0+\theta t_1,p}(\Omega))$ and we have~:
$$\parallel \Pi \parallel_{_{\mathcal{L}_{\theta}}}~\leq~
\parallel \Pi \parallel^{1-\theta}_{_{\mathcal{L}_0}}~\parallel \Pi\parallel^{\theta}_{_{\mathcal{L}_1}}~~~.$$
\end{theorem}

\noindent \textbf{Proof of Proposition \ref{proppit}} \\
As $H^{3/2}(\Gamma)$ is contained in $\mathcal{C}^0(\Gamma),~\Pi^{\Gamma}_{_{\mathcal{T}}}$ 
is well defined and, moreover, $\Pi^{\Gamma}_{_{\mathcal{T}}} \mu \in \Gamma_{_{\mathcal{T}}}$ if 
$\mu \in H^{3/2}(\Gamma)$. The difference   
$I - \Pi^{\Gamma}_{_{\mathcal{T}}}$ is continuous from $H^2(\Gamma)$ to $L^2(\Gamma)$ and its
norm is bounded by $C_1(\sigma)~h^2_{_{\mathcal{T}}}$ \cite{CR72}. It is also continuous 
from $H^1(\Gamma)$ to $H^1(\Gamma)$, its norm being bounded by a constant $C_2(\sigma)$. 
The interpolation theorem \ref{Iss} concludes
(with $\theta = 1/2$) that $I-\Pi^{\Gamma}_{_{\mathcal{T}}}$ is continuous from
$H^{3/2}(\Gamma)$ to $H^{1/2}(\Gamma)$ with its norm bounded by
$\sqrt{C_1(\sigma)h^2_{_{\mathcal{T}}}}\sqrt{C_2(\sigma)} = C'(\sigma)h_{_{\mathcal{T}}}$.$\hfill \blacksquare$

\begin{prop}{Regularity of components of functions in $M(\Omega)$.}\label{phi0h2} \\
Let us recall that the domain $\Omega$ is supposed to be convex. Let $\varphi$ be in 
$M(\Omega) \cap H^2(\Omega) = H^2(\Omega)$, splitted into : $\varphi~=~\varphi^0~+~\varphi^{\Delta}$
(see Proposition \ref{prop0}). Then $\varphi^0$ and $\varphi^{\Delta}$ belong to $H^2(\Omega)$.
\end{prop}

\proof \\
We recall that $\varphi^0 \in H^1_0(\Omega)$ is the solution of~:
$$
\left\{
\begin{array}{lll}
\Delta \varphi^0 & = & \Delta \varphi~\mbox{ in } \Omega \\
\gamma \varphi^0 & = & 0~~~~ \mbox{ on } \Gamma~~~.
\end{array}
\right.
$$
If $\varphi$ belongs to $H^2(\Omega)$, its laplacian belongs to $L^2(\Omega)$ and,  
by the regularity of the laplacian operator on a convex domain, $\varphi^0$ 
belongs to $H^2(\Omega)$. By difference, $\varphi^{\Delta} = \varphi - \varphi^0$ belongs also to $H^2(\Omega)$.
$\hfill \blacksquare$

\begin{defi}{Interpolation operators.} \label{def5}\\
The interpolation operator 
$\phi^{\infty}_{_{\mathcal{T}}} : \mathcal{H}(\Omega) \longrightarrow \mathcal{H}_{_{\mathcal{T},\infty}}$ 
is defined by $\phi^{\infty}_{_{\mathcal{T}}} \varphi^{\Delta} = \zeta$, where $\zeta \in H^1(\Omega)$ is the solution of
the Dirichlet problem~:
$$
\left\{
\begin{array}{llll}
\Delta \zeta & = & 0 & \mathrm{ in }~ \Omega \\
\gamma \zeta & = & \Pi^{\Gamma}_{_{\mathcal{T}}}(\gamma \varphi) & \mathrm{ on }~\Gamma~~~.
\end{array}
\right.
$$
We define the interpolation operator $\mathcal{P}_{_{\mathcal{T}}}$ from $M(\Omega) \cap H^2(\Omega) = H^2(\Omega)$ to 
$H^{1,\infty}_{\mathcal{T}} = H^1_{_{0,\mathcal{T}}} \oplus \mathcal{H}_{_{\mathcal{T},\infty}}$ by the relations~:
$$\mathcal{P}_{_{\mathcal{T}}} : H^2(\Omega) \ni \varphi = \varphi^0 + \varphi^{_{\Delta}} \longmapsto 
\mathcal{P}_{_{\mathcal{T}}} \varphi = \Pi_{_{\mathcal{T}}}\varphi^0 
+ \phi^{\infty}_{_{\mathcal{T}}} \varphi^{_{\Delta}} \in 
H^{1,\infty}_{_{\mathcal{T}}}~~~.$$ 
\end{defi}

\begin{prop}{Error estimates.} \label{prop5}\\
For $\mathcal{T} \in \mathcal{U}_{\sigma}$, we have the following estimates~:
\begin{equation} \label{phi0}
\parallel \varphi^0 - \Pi_{_{\mathcal{T}}}\varphi^0 \parallel_{_M}~\leq~
C(\sigma)~h_{_{\mathcal{T}}}~\parallel \varphi^0 \parallel_{_{2,\Omega}}
\end{equation}
\begin{equation} \label{phidelta}
\parallel \varphi^{\Delta} - \phi^{\infty}_{_{\mathcal{T}}} \varphi^{\Delta} \parallel_{_M}~\leq~
C(\sigma)~h_{_{\mathcal{T}}}~\parallel \varphi \parallel_{_{2,\Omega}}~~~.
\end{equation}
\end{prop}

\proof \\
\monitem As $H^2(\Omega) \subset \mathcal{C}^0(\overline{\Omega})$, we can use a classical interpolation operator 
and we have the following interpolation error results \cite{CR72}~:
$$\parallel \varphi^0 - \Pi_{_{\mathcal{T}}} \varphi^0 \parallel_{_{1,\Omega}}~\leq~
C~h_{_{\mathcal{T}}}~|\varphi^0|_{_{2,\Omega}}~~~.$$
But
$\parallel \varphi^0 - \Pi_{_{\mathcal{T}}} \varphi^0 \parallel_{_M}$ is equivalent to 
$\parallel \varphi^0 - \Pi_{_{\mathcal{T}}} \varphi^0 \parallel_{_{1,\Omega}}$ because
$\varphi^0 - \Pi_{_{\mathcal{T}}} \varphi^0$ belongs to $H^1_0(\Omega)$ (see Proposition \ref{rem2}). So~:
$$\parallel \varphi^0 - \Pi_{_{\mathcal{T}}} \varphi^0 \parallel_{_M}~\leq~
C~h_{_{\mathcal{T}}}~|\varphi^0|_{_{2,\Omega}}~~~,$$
and relation (\ref{phi0}) is established.\\
\monitem We now interpolate the harmonic part $\varphi^{\Delta}$ of $\varphi \in M(\Omega)$. 
As $\varphi$ is supposed to be in $H^2(\Omega)$, its trace $\gamma \varphi$ belongs to $H^{3/2}(\Gamma)$. 
Let $\lambda = \Pi^{\Gamma}_{_{\mathcal{T}}}(\gamma \varphi)$, 
$\lambda \in \Gamma_{\mathcal{T}}$ and we have (Proposition \ref{proppit})~:
\begin{equation}\label{pit2}
\parallel \gamma \varphi - \Pi^{\Gamma}_{_{\mathcal{T}}}(\gamma \varphi) \parallel_{_{1/2,\Gamma}}~=~
\parallel \gamma \varphi - \lambda \parallel_{_{1/2,\Gamma}}~\leq~
C~h_{_{\mathcal{T}}}~\parallel \gamma \varphi \parallel_{_{3/2,\Gamma}}~~~.
\end{equation}
Let $\zeta \in H^1(\Omega)$ be such that~:
\begin{eqnarray*}
\left\{
\begin{array}{lll}
\Delta \zeta & = & 0 \mbox{ in } \Omega \\
\gamma \zeta & = & \lambda \mbox{ on } \Gamma~~~.
\end{array}
\right.
\end{eqnarray*}
Thus, $C$ being various constants independent of $h_{_{\mathcal{T}}}$ but varying with $\sigma$, we obtain~:
\begin{eqnarray*}
\parallel \varphi^{\Delta}- \zeta \parallel_{_M} & = & \parallel \varphi^{\Delta} - \zeta
\parallel_{_{0,\Omega}} \mbox{ because } \varphi^{\Delta} \mbox{ and } \zeta \mbox{ are harmonic}\\
& \leq &  \parallel \varphi^{\Delta} - \zeta \parallel_{_{1,\Omega}} \\
& \leq & C \parallel \gamma \varphi^{\Delta} - \gamma \zeta \parallel_{_{1/2,\Gamma}} \mbox{ because trace 
controls harmonics}\\
& \leq & C \parallel \gamma \varphi - \lambda \parallel_{_{1/2,\Gamma}} \mbox{ as }
\gamma \varphi^{\Delta} = \gamma \varphi \mbox{ on } \Gamma \mbox{ and } \gamma \zeta = \lambda\\
& \leq & C h_{_{\mathcal{T}}} \parallel \gamma \varphi \parallel_{_{3/2,\Gamma}}
\mbox{ by relation (\ref{pit2})}\\
& \leq & C h_{_{\mathcal{T}}} \parallel \varphi \parallel_{_{2,\Omega}}
\mbox{ by continuity of the trace } \cite{LM}~~~,
\end{eqnarray*}
which gives the announced result. $\hfill \blacksquare$

\subsection{Error estimates}

We define an auxiliary problem~: 
\begin{eqnarray} \label{Pst}
\mbox{Find~} (\theta_{_{\mathcal{T}}},\eta_{_{\mathcal{T}}}) \in H^1_{_{0,\mathcal{T}}} 
\times H^{1,\infty}_{_{\mathcal{T}}} \mbox{~such~that~:} ~~~~~~~~~~~~~~~~~~~~~ \nonumber \\
\left\{
\begin{array}{lcll}
(\theta_{_{\mathcal{T}}},\varphi)~-~(\nabla \eta_{_{\mathcal{T}}},\nabla \varphi)
        & = & (g,\varphi)~+~\langle \Delta \varphi,m\rangle_{_{-1,1}} & \forall \varphi \in H^{1,\infty}_{_{\mathcal{T}}} \\
(\nabla \theta_{_{\mathcal{T}}},\nabla \xi)
        & = & \langle l,\xi\rangle_{_{-1,1}}                          & \forall \xi \in H^1_{_{0,\mathcal{T}}}~~~.
\end{array}
\right.
\end{eqnarray}

\begin{prop}{Stability of discrete formulation (\ref{E9-11}).} \label{prop4}\\
Let $g \in L^2(\Omega)$, $m \in H^1_0(\Omega)$, $l \in H^{-1}(\Omega)$ and 
$(\theta_{_{\mathcal{T}}},\eta_{_{\mathcal{T}}}) \in H^1_{_{0,\mathcal{T}}} \times H^{1,\infty}_{_{\mathcal{T}}}$ 
be the unique solution of problem (\ref{Pst}). We have stability in the following sense. 
There exists a constant $C$ only dependent on $\sigma$ such that the following stability inequality holds~:
$$\parallel \theta_{_{\mathcal{T}}} \parallel_{_M}~+~\parallel \nabla \eta_{_{\mathcal{T}}} \parallel_{_{0,\Omega}}
~\leq~C~\left\{\parallel g \parallel_{_{0,\Omega}}~+~\parallel \nabla m \parallel_{_{0,\Omega}}~+~
\parallel l \parallel_{_{-1,\Omega}}\right\}~~~.$$
\end{prop}

\proof \\
\monitem As in Proposition \ref{prop3}, three steps are required. We split each function of 
$H^{1,\infty}_{_{\mathcal{T}}}$ into two parts~:  
$\theta_{_{\mathcal{T}}} = \theta^0_{_{\mathcal{T}}} +\theta^{\Delta}_{_{\mathcal{T}}} \in 
H^1_{_{0,\mathcal{T}}} \oplus \mathcal{H}_{_{\mathcal{T},\infty}}$. And we obtain~:
\begin{eqnarray}
(\theta^0_{_{\mathcal{T}}},\varphi^0) +
(\theta^{\Delta}_{_{\mathcal{T}}},\varphi^0)
- (\nabla \eta_{_{\mathcal{T}}},\nabla \varphi^0) \! \! \! \! & = & \! \! \! \!
(g,\varphi^0) - \!\langle \! \Delta \varphi^0,m \!\rangle_{_{-1,1}} \! \forall \varphi^0 \! \in \! H^1_{_{0,\mathcal{T}}}
\label{E33a} \\
(\theta^0_{_{\mathcal{T}}},\varphi^{\Delta}) +
(\theta^{\Delta}_{_{\mathcal{T}}},\varphi^{\Delta})
\! \! \! \! & = & \! \! \! \! (g,\varphi^{\Delta})~~~~~~~\forall \varphi^{\Delta} \in \mathcal{H}_{_{\mathcal{T},\infty}} \label{E33b} \\
(\nabla \theta^0_{_{\mathcal{T}}}, \nabla \xi) \! \! \! \! & = & \! \! \! \! \langle l,\xi\rangle_{_{-1,1}}~\forall \xi
\in H^1_{_{0,\mathcal{T}}}~~~. \label{E33c}
\end{eqnarray}
As $m \in H^1_0(\Omega)$), notice that~:
$$- \langle \Delta \varphi^0,m\rangle_{_{-1,1}}~=~(\nabla \varphi^0, \nabla m)~~,~~\forall \varphi^0 \in H^1_{_{0,\mathcal{T}}}~~~.$$
In the following, $C$ will denote various constants. \\
\monitem {\bf{First step.}} In the last equation (\ref{E33c}), we take $\xi = \theta^0_{_{\mathcal{T}}}$ 
and, as $\theta^0_{_{\mathcal{T}}}$ belongs to $H^1_0(\Omega)$, we obtain~:
$$\parallel \nabla \theta^0_{_{\mathcal{T}}} \parallel_{_{0,\Omega}}~\leq~C~\parallel l \parallel_{_{-1,\Omega}}~~~,$$
because of the Poincar\'e's inequality. Then, the equivalence of $H^1_0-$norm and $M-$norm for functions in 
$H^1_0(\Omega)$ leads to~:
$$\parallel \theta^0_{_{\mathcal{T}}} \parallel_{_M}~\leq~C~\parallel l \parallel_{_{-1,\Omega}}~~~.$$
\monitem {\bf{Second step.}} In the equation (\ref{E33b}), we take $\varphi^{\Delta} =
\theta^{\Delta}_{_{\mathcal{T}}}$. As $\Delta \theta^{\Delta}_{_{\mathcal{T}}} = 0$ and 
$\theta^0_{_{\mathcal{T}}}$ is fixed by the previous step, we obtain~: 
$$\parallel \theta^{\Delta}_{_{\mathcal{T}}} \parallel_{_M}
~=~\parallel \theta^{\Delta}_{_{\mathcal{T}}} \parallel_{_{0,\Omega}}
~\leq~\parallel g \parallel_{_{0,\Omega}}~+~\parallel \theta^0_{_{\mathcal{T}}} \parallel_{_{0,\Omega}}
~\leq~C~(\parallel g \parallel_{_{0,\Omega}}~+~\parallel l \parallel_{_{-1,\Omega}})$$
\monitem {\bf{Third step.}} Finally, in the first equation (\ref{E33a}), $\theta^0_{_{\mathcal{T}}}$ and 
$\theta^{\Delta}_{_{\mathcal{T}}}$ being known, we take $\varphi^0 = \eta_{_{\mathcal{T}}}$ and obtain~:
$$\parallel \nabla \eta_{_{\mathcal{T}}} \parallel_{_{0,\Omega}}~\leq~
C~\left( \parallel \theta^0_{_{\mathcal{T}}} \parallel_{_{0,\Omega}}~+~
\parallel \theta^{\Delta}_{_{\mathcal{T}}} \parallel_{_{0,\Omega}}~+~\parallel g \parallel_{_{0,\Omega}} \right)
~+~\parallel \nabla m \parallel_{_{0,\Omega}}~~~.$$
By combining the three steps, Proposition \ref{prop4} is proved. $\hfill \blacksquare$

\begin{prop}{Interpolation error majorates the error.} \label{ieme} \\
Let $(\omega,\psi)$ be the solution of the continuous problem (\ref{P}) and 
$(\omega_{_{\mathcal{T}}},\psi_{_{\mathcal{T}}})$ the solution of (\ref{E9-11}). We have~:  
$$\parallel \omega - \omega_{_{\mathcal{T}}}  \parallel_{_M}~+~
\parallel \psi - \psi_{_{\mathcal{T}}} \parallel_{_{1,\Omega}}
~\leq~C~( \parallel \omega - \mathcal{P}_{_{\mathcal{T}}} \omega \parallel_{_M}
~+~\parallel \psi - \Pi_{_{\mathcal{T}}} \psi \parallel_{_{1,\Omega}})~~~.$$
\end{prop}

\proof \\
\monitem The continuous problem (\ref{P}) can be written with a test function $\varphi$ in $H^{1,\infty}_{_{\mathcal{T}}}$
since $H^{1,\infty}_{_{\mathcal{T}}} \subset M(\Omega)$ and with $\xi \in H^1_{_{0,\mathcal{T}}} \subset
H^1_0(\Omega)$~:
\begin{equation} \label{Pt}
\left\{
\begin{array}{llll}
(\omega,\varphi)~+~\langle \Delta \varphi,\psi\rangle_{_{-1,1}} & = & 0 & \forall \varphi \in H^{1,\infty}_{_{\mathcal{T}}}\\
\langle -\Delta \omega,\xi\rangle_{_{-1,1}} & = & (\mathbf{f},\mathbf{rot}~\xi) & \forall \xi \in H^1_{_{0,\mathcal{T}}}
\end{array}
\right.
\end{equation}
We subtract (\ref{Pt}) and (\ref{E9-11})~:
\begin{eqnarray}\label{diff}
\left\{
\begin{array}{llll}
(\omega - \omega_{_{\mathcal{T}}},\varphi)~+~
\langle \Delta \varphi,\psi\rangle_{_{-1,1}}~+~(\nabla \psi_{_{\mathcal{T}}},\nabla \varphi) & = & 0 
& \forall \varphi \in H^{1,\infty}_{_{\mathcal{T}}}\\
\langle -\Delta \omega,\xi\rangle_{_{-1,1}}~-~(\nabla \omega_{_{\mathcal{T}}},\nabla \xi)    & = & 0 
& \forall \xi \in H^1_{_{0,\mathcal{T}}}
\end{array}
\right.
\end{eqnarray}
We introduce $\Pi_{_{\mathcal{T}}} \psi$ that interpolates $\psi$ on the mesh $\mathcal{T}$
and $\mathcal{P}_{_{\mathcal{T}}} \omega$ that interpolates
$\omega \in M(\Omega)$ on $\mathcal{T}$ (these interpolants, both in $H^1(\Omega)$, are defined in 
Proposition \ref{prop5}). \\
\monitem Let us begin with the first equation of problem (\ref{diff}). We add and subtract 
$\mathcal{P}_{_{\mathcal{T}}} \omega$ and $\Pi_{_{\mathcal{T}}} \psi$ and obtain, for all $\varphi$ 
in $H^{1,\infty}_{_{\mathcal{T}}}$~: 
$$(\omega - \mathcal{P}_{_{\mathcal{T}}} \omega,\varphi) 
~+~(\mathcal{P}_{_{\mathcal{T}}} \omega - \omega_{_{\mathcal{T}}},\varphi) 
~+~\langle \Delta \varphi,\psi - \Pi_{_{\mathcal{T}}} \psi \rangle_{_{-1,1}}~~~~~~~~~~~~~~~~~~~~~~$$
$$~+~\langle \Delta \varphi,\Pi_{_{\mathcal{T}}} \psi \rangle_{_{-1,1}}
~+~(\nabla (\psi_{_{\mathcal{T}}} - \Pi_{_{\mathcal{T}}} \psi),\nabla \varphi)
~+~(\nabla \Pi_{_{\mathcal{T}}} \psi,\nabla \varphi)~=~0~~~,$$
or else~:
$$(\mathcal{P}_{_{\mathcal{T}}} \omega - \omega_{_{\mathcal{T}}},\varphi)
~-~(\nabla (\Pi_{_{\mathcal{T}}} \psi - \psi_{_{\mathcal{T}}}),\nabla \varphi)
~=~-(\omega - \mathcal{P}_{_{\mathcal{T}}} \omega,\varphi)~~~~~~~~~~~~~~~~~~$$
$$- \langle \Delta \varphi,\psi - \Pi_{_{\mathcal{T}}} \psi \rangle_{_{-1,1}} 
\underbrace{- \langle \Delta \varphi,\Pi_{_{\mathcal{T}}} \psi \rangle_{_{-1,1}} 
- (\nabla \Pi_{_{\mathcal{T}}} \psi,\nabla \varphi)}_{= 0}~~~.$$
Finally, the first equation becomes~:
\begin{eqnarray} \label{diff1} 
(\mathcal{P}_{_{\mathcal{T}}} \omega - \omega_{_{\mathcal{T}}},\varphi)
& - & (\nabla (\Pi_{_{\mathcal{T}}} \psi - \psi_{_{\mathcal{T}}}),\nabla \varphi) \nonumber \\
& = & -(\omega - \mathcal{P}_{_{\mathcal{T}}} \omega,\varphi) 
~-~\langle \Delta \varphi,\psi - \Pi_{_{\mathcal{T}}} \psi \rangle_{_{-1,1}}~~~.
\end{eqnarray} 
\monitem Using the same techniques for the second equation, we obtain, for all $\xi$ 
in $H^1_{_{0,\mathcal{T}}}$~: 
$$\langle -\Delta (\omega - \mathcal{P}_{_{\mathcal{T}}} \omega),\xi\rangle_{_{-1,1}}~+~
\langle -\Delta \mathcal{P}_{_{\mathcal{T}}} \omega ,\xi\rangle_{_{-1,1}}~~~~~~~~~~~~~~~~~~$$
$$-(\nabla (\omega_{_{\mathcal{T}}}~-~\mathcal{P}_{_{\mathcal{T}}} \omega),\nabla \xi)
~-~(\nabla \mathcal{P}_{_{\mathcal{T}}} \omega,\nabla \xi)~=~0~~~,$$
and then~:
\begin{eqnarray}
\begin{array}{rcl}\label{diff2} 
(\nabla (\mathcal{P}_{_{\mathcal{T}}} \omega - \omega_{_{\mathcal{T}}}),\nabla \xi)
& = & \underbrace{(\nabla \mathcal{P}_{_{\mathcal{T}}} \omega,\nabla \xi) -
\langle -\Delta \mathcal{P}_{_{\mathcal{T}}} \omega ,\xi\rangle_{_{-1,1}}}_{= 0} \\
& - & \langle -\Delta (\omega - \mathcal{P}_{_{\mathcal{T}}} \omega),\xi\rangle_{_{-1,1}}~~~.
\end{array}
\end{eqnarray}
\monitem Equations (\ref{diff1}) and (\ref{diff2}) lead to the following problem~:
\begin{eqnarray*}
\left\{
\begin{array}{rcl}
(\mathcal{P}_{_{\mathcal{T}}} \omega - \omega_{_{\mathcal{T}}},\varphi)
& - & (\nabla (\Pi_{_{\mathcal{T}}} \psi - \psi_{_{\mathcal{T}}}),\nabla \varphi)\\
& = & -(\omega - \mathcal{P}_{_{\mathcal{T}}} \omega,\varphi)
~-~ \! \langle \Delta \varphi,\psi - \Pi_{_{\mathcal{T}}} \psi \rangle_{_{-1,1}}
~,~\forall \varphi \in H^{1,\infty}_{_{\mathcal{T}}} \\
(\nabla (\mathcal{P}_{_{\mathcal{T}}} \omega - \omega_{_{\mathcal{T}}}),\nabla \xi)
& = & \langle \Delta (\omega - \mathcal{P}_{_{\mathcal{T}}} \omega),\xi\rangle_{_{-1,1}} 
~,~\forall \xi \in H^1_{_{0,\mathcal{T}}}
\end{array}
\right.
\end{eqnarray*}
It is the auxiliary problem (\ref{Pst}) with~: 
$g = \mathcal{P}_{_{\mathcal{T}}} \omega - \omega \in L^2(\Omega)$, $m = \Pi_{_{\mathcal{T}}} \psi -
\psi \in H^1_0(\Omega)$, $l =\Delta (\omega - \mathcal{P}_{_{\mathcal{T}}} \omega) \in H^{-1}(\Omega)$. 
The triangular inequality and Proposition \ref{prop4} lead to~:
\begin{eqnarray*}
\parallel \omega - \omega_{_{\mathcal{T}}}  \parallel_{_M} +
\parallel \psi - \psi_{_{\mathcal{T}}} \parallel_{_{1,\Omega}} 
& \leq & \underbrace{\parallel \omega - \mathcal{P}_{_{\mathcal{T}}} \omega \parallel_{_M} +
\parallel \psi - \Pi_{_{\mathcal{T}}} \psi \parallel_{_{1,\Omega}}}_{\mathrm{interpolation~error}} \\
& + & \parallel \mathcal{P}_{_{\mathcal{T}}} \omega - \omega_{_{\mathcal{T}}} \parallel_{_M}
+ \parallel \Pi_{_{\mathcal{T}}} \psi - \psi_{_{\mathcal{T}}} \parallel_{_{1,\Omega}} 
\end{eqnarray*}
\begin{eqnarray*}
& \leq & \parallel \omega - \mathcal{P}_{_{\mathcal{T}}} \omega \parallel_{_M} 
         ~+~\parallel \psi - \Pi_{_{\mathcal{T}}} \psi \parallel_{_{1,\Omega}} \\
& ~    & +~C~\left( \parallel \mathcal{P}_{_{\mathcal{T}}} \omega - \omega \parallel_{_{0,\Omega}}
          + \parallel \Delta (\omega - \mathcal{P}_{_{\mathcal{T}}} \omega) \parallel_{_{-1,\Omega}}
	  + \parallel \nabla(\psi - \Pi_{_{\mathcal{T}}} \psi) \parallel_{_{0,\Omega}} \right)
\end{eqnarray*}
and finally~: 
\begin{eqnarray*}
\parallel \omega - \omega_{_{\mathcal{T}}}  \parallel_{_M}~+~
\parallel \psi - \psi_{_{\mathcal{T}}} \parallel_{_{1,\Omega}}
~\leq~C~( \parallel \omega - \mathcal{P}_{_{\mathcal{T}}} \omega \parallel_{_M}
~+~\parallel \psi - \Pi_{_{\mathcal{T}}} \psi \parallel_{_{1,\Omega}})
\end{eqnarray*}
which achieves the proof. $\hfill \blacksquare$\\

\noindent The previous Proposition and the interpolation error lead to the following result.

\begin{theorem}{Convergence result.} \label{th1} \\
If $\mbox{rot}~\mathbf{f}$ belongs to $L^2(\Omega)$ and $\mathcal{T}$ belongs to a regular family of 
triangulation $\mathcal{U}_{\sigma}$, the discrete problem (\ref{E9-11}) has a unique solution 
$(\psi_{_{\mathcal{T}}},\omega_{_{\mathcal{T}}})$ in $H^1_{_{0,\mathcal{T}}} \times H^{1,\infty}_{_{\mathcal{T}}}$,
associated with a stable discretization of the Stokes problem (\ref{E1})-(\ref{E4}). 
If $(\psi,\omega) \in H^2(\Omega) \cap H^1_0(\Omega) \times H^2(\Omega)$, there exists a strictly 
positive constant $C(\sigma)$ such that, for all $\mathcal{T}$ in $\mathcal{U}_{\sigma}$, we have~:
$$\parallel \omega - \omega_{_{\mathcal{T}}}  \parallel_{_M}~+~\parallel \psi - \psi_{_{\mathcal{T}}} \parallel_{_{1,\Omega}}
~\leq~C(\sigma)~h_{_{\mathcal{T}}}~(\parallel \omega \parallel_{_{2,\Omega}} + \parallel \psi \parallel_{_{2,\Omega}})~~~.$$
\end{theorem}

\proof \\
Let us recall that, due to Proposition \ref{prop0}, $\omega^0$ is such that~:
$$
\left\{
\begin{array}{lll}
\Delta \omega^0 & = & \Delta \omega~=~\mbox{rot}~\mathbf{f} \mbox{ in } \Omega \\
\gamma \omega^0 & = & 0~~~ \mbox{ on } \Gamma~~~.
\end{array}
\right.
$$
As $\Omega$ is a convex polygon, and as $\mbox{rot}~\mathbf{f}$ belongs to
$L^2(\Omega)$, $\omega^0$ belongs to $H^2(\Omega)$ \cite{Gri} and we have~: 
\begin{equation} \label{maj}
\parallel \omega^0 \parallel_{_{2,\Omega}}~\leq~C~\parallel \Delta \omega \parallel_{_{0,\Omega}} 
~\leq~C~\parallel \omega \parallel_{_{2,\Omega}}~~~.
\end{equation}
So, by the help of Proposition (\ref{ieme}), $C$ denoting various constants, we obtain~: 
\begin{eqnarray*}
\parallel \omega - \omega_{_{\mathcal{T}}} \parallel_{_M} 
& + & \parallel \psi - \psi_{_{\mathcal{T}}} \parallel_{_{1,\Omega}}~\leq~
      C~\left( \parallel \omega - \mathcal{P}_{_{\mathcal{T}}} \omega \parallel_{_M} + 
      \parallel \psi - \Pi_{_{\mathcal{T}}} \psi \parallel_{_{1,\Omega}} \right) \\
& \leq & C \left( \parallel \omega^0 + \omega^{\Delta} - \Pi_{_{\mathcal{T}}} \omega^0 
- \phi^{\infty}_{_{\mathcal{T}}} \omega^{\Delta} \parallel_{_M} 
+ h_{_{\mathcal{T}}} \parallel \psi \parallel_{_{2,\Omega}} \right)\\
& \leq & C \left( \parallel \omega^0 - \Pi_{_{\mathcal{T}}} \omega^0 \parallel_{_M}
+ \parallel \omega^{\Delta} - \phi^{\infty}_{_{\mathcal{T}}} \omega^{\Delta} \parallel_{_M} 
+ h_{_{\mathcal{T}}} \parallel \psi \parallel_{_{2,\Omega}} \right) \\
& \leq & C h_{_{\mathcal{T}}} \left( \parallel \omega^0 \parallel_{_{2,\Omega}} 
+ \parallel \omega \parallel_{_{2,\Omega}}
+ \parallel \psi \parallel_{_{2,\Omega}} \right) \mbox{ by Proposition \ref{prop5}}\\
& \leq & C h_{_{\mathcal{T}}} \left( \parallel \omega \parallel_{_{2,\Omega}} 
+ \parallel \omega \parallel_{_{2,\Omega}} + 
\parallel \psi \parallel_{_{2,\Omega}} \right) \mbox{ by (\ref{maj})} \\
& \leq & C h_{_{\mathcal{T}}} \left( \parallel \omega \parallel_{_{2,\Omega}} 
+ \parallel \psi \parallel_{_{2,\Omega}} \right)~~~, 
\end{eqnarray*}
which gives the announced result.$\hfill \blacksquare$\\

\begin{rem}{~~~}\\
\monitem The theorem \ref{th1} is important because it shows that using a space of harmonic functions
along the boundary gives an error of order $\mathcal{O}(h_{_{\mathcal{T}}})$ when $\omega \in H^2(\Omega)$. 
It is now necessary to approach this space numerically {\bf{without losing optimality}}, 
\emph{i.e} error between discrete harmonic functions and real harmonic functions must be
bounded by $h_{_{\mathcal{T}}}$. We propose in this paper a way of approaching 
$\mathcal{H}_{_{\mathcal{T},\infty}}$ by use of mesh refinements but other possibilities are in progress 
(\cite{ASS}, \cite{Ru1}).\\
\monitem It seems possible to reduce hypothesis on the regularity of the vorticity and to obtain 
convergence of order $\mathcal{O}(h_{_{\mathcal{T}}})$ for $\omega$ only in $H^1(\Omega)$.
\end{rem}

\newpage 
\section{Discrete harmonic functions}

\subsection{Formulation}

Let $\mathcal{T} \in \mathcal{U}_{\sigma}$ be a given mesh. Then, the family of meshes 
$\mathcal{T}_k$ ($k \in \mbox{I\hspace{-.15em}N}$) is defined recursively as follows~:
\begin{itemize}
\item[(i)] $\mathcal{T}_0 = \mathcal{T}$;
\item[(ii)] $\mathcal{T}_{k+1} \subset \mathcal{T}_k$ is obtained from $\mathcal{T}_k$ by
dividing each triangle in four homothetic ones (see Figure~\ref{recursive}).
\end{itemize}
%
\begin{figure}    [H]  \centering
\centerline  {\includegraphics[width=.30\textwidth]   {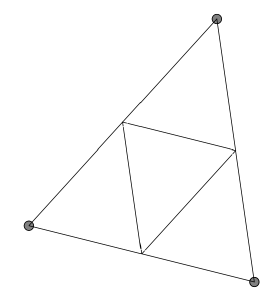}}
\caption{Division of a triangle of $\mathcal{T}_k$.}
\label{recursive}
\end{figure}
%
It is clear that $(\mathcal{T}_k)_{k \in \N}$ belong to $\mathcal{U}_{\sigma}$ defined above. 
We denote $H^1_{_{\mathcal{T}_k}}$ the space of all 
continuous functions on $\overline{\Omega}$, polynomial of degree $1$ in each triangle of
$\mathcal{T}_{k}$ and $H^1_{_{0,\mathcal{T}_k}} = H^1_{_{\mathcal{T}_k}} \cap H^1_0(\Omega)$.

\begin{defi}{Discrete harmonics.}\\
We define the space $\mathcal{H}_{_{\mathcal{T},k}}$ of discrete harmonic functions by~:
$$\mathcal{H}_{_{\mathcal{T},k}} = \left\{ \varphi \in
H^1_{_{\mathcal{T}_k}},
\exists \lambda \in \Gamma_{_{\mathcal{T}}}~~\gamma \varphi =
\lambda \mathrm{~on~} \Gamma,(\nabla \varphi,\nabla \xi) = 0~~
\forall \xi \in H^1_{_{0,\mathcal{T}_k}} \right\}~~~.$$
\end{defi}
We define the following lifting operator $Z_k$ of order $k$~: 
$$Z_k : \Gamma_{_{\mathcal{T}}} \ni \lambda \longmapsto Z_k(\lambda) \in \mathcal{H}_{_{\mathcal{T},k}}~~~,$$
by the relations~: 
\begin{eqnarray*}
\left\{
\begin{array}{rcl}
(\nabla Z_k(\lambda),\nabla \xi) & = & 0~~\forall \xi \in H^1_{_{0,\mathcal{T}_k}}\\
\gamma Z_k(\lambda) & = & \lambda~~\mbox{ on } \Gamma~~~.
\end{array}
\right.
\end{eqnarray*}
The space $\mathcal{H}_{_{\mathcal{T},k}}$ is finite-dimensional and $\mbox{dim }
\mathcal{H}_{_{\mathcal{T},k}} = \mbox{dim } \Gamma_{_{\mathcal{T}}}$. The nodal values on 
$\mathcal{T}_k \cap \Gamma$ are obtained by linear interpolation from values at 
the vertices of the initial mesh $\mathcal{T}_0 = \mathcal{T}.$ Internal values of function 
$Z_k(\lambda)$ are obtained by solving a discrete Laplace equation in space $H^1_{_{0,\mathcal{T}_k}}$ 
with boundary values from the previous step.\\ 

\vfill \eject ~ \smallskip
\begin{figure}    [H]  \centering
\centerline  {\includegraphics[width=.50\textwidth]   {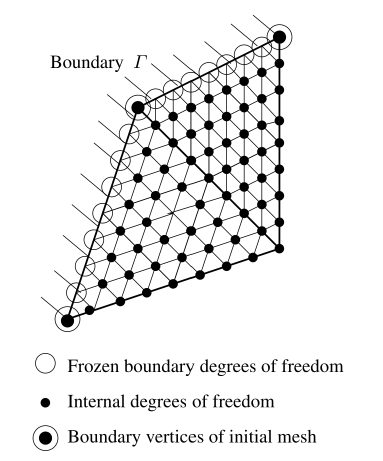}}
\caption{Mesh for discrete harmonics at level $k = 3$.}
\label{ddlrec}
\end{figure}
%
\begin{prop}{Continuous decomposition of space $\mathcal{H}_{_{\mathcal{T},k}}$.} \label{prop6}\\
For all $\lambda \in \Gamma_{_{\mathcal{T}}}$,
we have the following decomposition~:
$$Z_k(\lambda) = Z_k(\lambda)^0 + Z_{\infty}(\lambda)~~~,$$
where $Z_{\infty}(\lambda)$, which is independent of $k$, is defined in Section $3.1$ 
and $Z_k(\lambda)^0$ belongs to $H^1_0(\Omega)$. Moreover, we have~:
\begin{equation} \label{norme-1}
\parallel \Delta Z_k(\lambda) \parallel_{_{-1,\Omega}}~=~
\parallel \nabla Z_k(\lambda)^0 \parallel_{_{0,\Omega}}~~~.
\end{equation}
\end{prop}

\proof \\
\monitem Let $Z_{\infty}(\lambda) \in \mathcal{H}_{_{\mathcal{T},\infty}}$ be such that
$\gamma Z_{\infty}(\lambda) = \lambda$ on $\Gamma$ (see problem (\ref{zinfini})). So, 
$\gamma (Z_k(\lambda) - Z_{\infty}(\lambda)) = 0$ on $\Gamma$. Then, 
the field $Z_k(\lambda)^0  = Z_k(\lambda) - Z_{\infty}(\lambda)$ belongs to $H^1_0(\Omega)$.\\
\monitem By definition, we have~:
$$\parallel \Delta Z_k(\lambda) \parallel_{_{-1,\Omega}}~=~
\sup_{\xi \in H^1_{_0}(\Omega) \backslash \left\{ 0 \right\}} 
\frac{\langle \Delta Z_k(\lambda),\xi \rangle_{_{-1,1}}}{\parallel \nabla \xi \parallel_{_{0,\Omega}}}~~~.$$
As $Z_k(\lambda)$ belongs to $H^1(\Omega)$, we may rewrite the duality product~:
$$\parallel \Delta Z_k(\lambda) \parallel_{_{-1,\Omega}}~=~
\sup_{\xi \in H^1_{_0}(\Omega)\backslash \left\{ 0 \right\}}
\frac{(\nabla Z_k(\lambda),\nabla \xi)}{\parallel \nabla \xi \parallel_{_{0,\Omega}}}
~=~\sup_{\xi \in H^1_{_0}(\Omega)\backslash \left\{ 0 \right\}}
\frac{(\nabla Z_k(\lambda)^0 ,\nabla \xi)}{\parallel \nabla \xi \parallel_{_{0,\Omega}}}~~~.$$
Finally, as $Z_k(\lambda)^0$ belongs to $H^1_0(\Omega)$, we can choose $\xi = Z_k(\lambda)^0$
in the above relation and the equality (\ref{norme-1}) is established.$\hfill \blacksquare$\\

To illustrate the non-harmonic part $Z_k(\lambda)^0$ of functions in 
$\mathcal{H}_{_{\mathcal{T},k}}$, we made some numerical experiments.
Let $\mathcal{T}$ be the mesh named $A$ described in Figure \ref{f:meshab}, $S$ be a vertex of 
$\mathcal{S}_{_{\mathcal{T}}} \cap \Gamma$ and $h_{_S} \in \mathcal{H}_{_{\mathcal{T},0}}$ 
be a discrete harmonic function of level $0$, used \emph{e.g.} in Glowinski-Pironneau \cite{GP}, and defined by~:  
\begin{eqnarray*}
\left\{
\begin{array}{rcll}
(\nabla h_{_S},\nabla\xi) & = & 0 & \forall \xi \in H^1_{_{0,\mathcal{T}}} \\
h_{_S}(T) & = & \delta_{_{ST}} & \mbox{ on } \Gamma~~~.
\end{array}
\right.
\end{eqnarray*}
It means $\gamma h_{_S} = \lambda_{S}$ where $\lambda_{S}$ is the corresponding basis function 
(see Figure \ref{lambda}). The previous proposition gives the following decomposition~:
$$h_{_S} =  Z_0(\lambda_{S})^0 + Z_{\infty}(\lambda_{S})~~~.$$
Let $\eta_{S,k} \in H^1_{_{0,\mathcal{T}_k}}$ be an approximation of the non harmonic 
part $Z_0(\lambda_{S})^0$ of $h_{_S}$, defined on the mesh $\mathcal{T}_k$ by~: 
\begin{eqnarray*}
\left\{
\begin{array}{rcll}
(\nabla \eta_{S,k},\nabla \chi) & = & (\nabla h_{_S},\nabla \chi) & \forall \chi \in H^1_{_{0,\mathcal{T}_k}} \\
\eta_{S,k} & = & 0 & \mbox{ on } \Gamma~~~.
\end{array}
\right.
\end{eqnarray*}
As $Z_0(\lambda_{S})^0$ verifies~:
$$(\nabla Z_0(\lambda_{S})^0,\nabla \chi)~=~(\nabla h_{_S},\nabla \chi)~~,~~\forall \chi \in H^1_0(\Omega)~~~,$$
we have~: $\eta_{S,k} \stackrel{k \longrightarrow \infty}{\longrightarrow} Z_0(\lambda_{S})^0$.\\
The complete shape of $\eta_{S,k}$ is well described, as we will prove further,  
after the {\bf{third}} refinement ($k=3$). 
The next figure (Figure \ref{f:fhd}) is a vizualisation of the function $\eta_{S,3}$ on the 
initial mesh $\mathcal{T} = \mathcal{T}_0$. We observe that all the "structure" of the function 
(minimum value $= -0.12$, maximum value $=0.19$) is confined inside the three triangles containing the 
vertex~$S$. 

\begin{figure}    [H]  \centering
\centerline  {\includegraphics[width=.90\textwidth]   {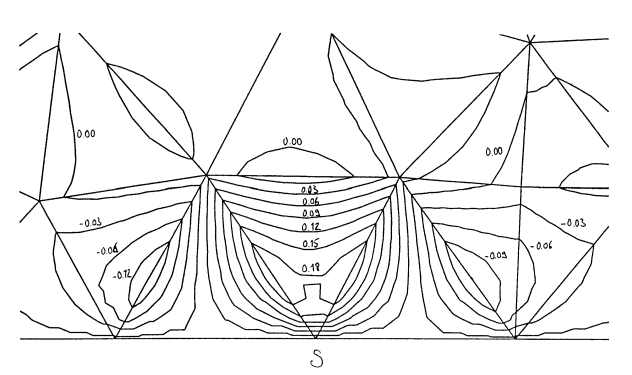}}
\caption{Vizualisation of $\eta_{S,3}$ (smiling function)}
\label{f:fhd}
\end{figure}

\smallskip

\begin{prop}{Interpolation error in the family of recursive meshes.} \label{prop9}\\
If $\mathcal{T}$ belongs to a regular family of triangulation $\mathcal{U}_{\sigma}$,
let $\Pi_{_{\mathcal{T}}} \varphi$ be the Lagrange interpolate on $\mathcal{T}$, which is well 
defined on $H^{1+\eta}(\Omega)$ for $\eta > 0$. Then, we have~: 
\begin{eqnarray} \label{inttk}
\forall~\varepsilon~\in~]0,1[~,~\forall \eta~\in~]0,1-\varepsilon[~,~
\forall~\varphi~\in~H^{2-\varepsilon}(\Omega)~,~\exists~\theta~\in~[\varepsilon,1]~, \nonumber \\
\exists~C(\sigma,\eta,\varepsilon)~>~0 \mbox{ such that }
\parallel \varphi -  \Pi_{_{\mathcal{T}}} \varphi \parallel_{_{1,\Omega}}
~\leq~C~h^{1-\theta}_{_{\mathcal{T}}}~\parallel \varphi \parallel_{_{2-\varepsilon,\Omega}}
\end{eqnarray}
\end{prop}

\begin{rem}
This proposition looks technical because the needed "classical" interpolation result~: 
$\parallel \varphi -  \Pi_{_{\mathcal{T}}} \varphi \parallel_{_{1,\Omega}}
\leq C h_{_{\mathcal{T}}} \parallel \varphi \parallel_{_{2,\Omega}}$ is not correct as  
functions in $\mathcal{H}_{_{\mathcal{T},\infty}}$ are not in $H^2(\Omega)$ but 
only in $H^{2-\varepsilon}(\Omega)$ (see Proposition \ref{prop8}). Nevertheless, it would
be possible to use the Scott and Zhang regularization for which functions only need to be 
in $H^1(\Omega)$ (see \cite{SZ}).
\end{rem}

\proof \\
Because $H^{2-\varepsilon}(\Omega) \subset \mathcal{C}^0(\overline{\Omega})$, it is possible 
to define the classical interpolate of $\varphi$ on the mesh $\mathcal{T}$. 
But $\varphi$ belongs only to $H^{2-\varepsilon}(\Omega)$. 
For the real number $2-\varepsilon$, we can find a parameter $\theta$ such that~:
 $2-\varepsilon = \theta (1+\eta) + 2 (1-\theta)$. Indeed $\theta$ 
is such that $0 < \theta < 1$, if $1+\eta < 2-\varepsilon$, and is defined by~:
\begin{equation}\label{theta}
\displaystyle{\theta = \frac{\varepsilon}{1-\eta}}~~~.
\end{equation}
Then, we have the following imbeddings~: $H^1(\Omega) \subset H^{1+\eta}(\Omega) \subset H^{2-\varepsilon}(\Omega) 
\subset \mathcal{C}^0(\overline{\Omega})$. The difference $I-\Pi_{_{\mathcal{T}}}$ is continuous 
from $H^2(\Omega)$ to $H^1(\Omega)$ with its norm bounded by $C_1(\sigma)h_{_{\mathcal{T}}}$.  
It is also continuous from $H^{1+\eta}(\Omega)$ to $H^1(\Omega)$
with its norm bounded by a constant $C_2(\sigma,\eta)$. So, for $\theta$ chosen as in (\ref{theta}), 
$I-\Pi_{_{\mathcal{T}}}$ is continuous from $H^{2-\varepsilon}(\Omega)$ to $H^1(\Omega)$ and its norm 
is bounded by $\displaystyle{(C_1(\sigma)h_{_{\mathcal{T}}})^{1-\theta}
  (C_2(\sigma,\eta))^{\theta} =}  $ 
$\displaystyle{(  = C(\sigma,\eta,\varepsilon) h^{1-\theta}_{_{\mathcal{T}}}}$, which leads to (\ref{inttk}).
$\hfill \blacksquare$

\begin{theorem} {Stability of discrete harmonics.} \label{th2}\\
Let $\delta$ be in $]0;+\infty[$. For any mesh $\mathcal{T}$ of the family $\mathcal{U}_{\sigma}$, there
exists an integer $K(\delta,\mathcal{T})$, such that~:
\begin{equation}\label{delta} 
\forall k \geq K(\delta,\mathcal{T})~,~\forall \varphi \in \mathcal{H}_{_{\mathcal{T},k}}~:~
\parallel \varphi \parallel_{_{0,\Omega}} \, \geq \, \delta \parallel \Delta
\varphi \parallel_{_{-1,\Omega}}~~~.
\end{equation}
\end{theorem}

\begin{rem}
Theorem \ref{th2} means that if $k$ is large enough ($k \geq 3$ as observed previously with the 
structure of functions $\eta_{S,3}$), the $L^2$-norm is sufficient to control the
$M$-norm of the vorticity. In other words, the laplacian of 
discrete harmonic functions, which are not rigourously harmonic, goes to zero as $k$ becomes large enough.
\end{rem}

\proof \\
We study the ratio $\displaystyle{\frac{\parallel \Delta \varphi \parallel_{_{-1,\Omega}}}
{\parallel \varphi \parallel_{_{0,\Omega}}}}$ for $\varphi$ arbitrary in 
$\mathcal{H}_{_{\mathcal{T},k}}$. First, we can write $\varphi = Z_k(\lambda)$ for an arbitrary discrete 
boundary function $\lambda \in \Gamma_{_{\mathcal{T}}}$ and we set 
$\displaystyle{\rho(\lambda) = \frac{\parallel \Delta Z_k(\lambda) \parallel_{_{-1,\Omega}}}
{\parallel Z_k(\lambda) \parallel_{_{0,\Omega}}}}$. \\
\monitem We first bound the numerator of $\rho(\lambda)$~: 
\begin{eqnarray*}
\parallel \Delta Z_k(\lambda) \parallel_{_{-1,\Omega}} 
& =    & \parallel \nabla Z_k(\lambda)^0 \parallel_{_{0,\Omega}}  \mbox{ by relation } (\ref{norme-1})\\
& \leq & \frac{1}{\sqrt{1+C^2_{p}}} \parallel Z_k(\lambda)^0 \parallel_{_{1,\Omega}} 
          (C_p \mbox{ is the Poincar\'e's constant}) \\
& \leq & \frac{1}{\sqrt{1+C^2_{p}}} \parallel Z_{\infty}(\lambda) - Z_k(\lambda) \parallel_{_{1,\Omega}}~~~.
\end{eqnarray*}
We have used~: $Z_k(\lambda) = Z_k(\lambda)^0 + Z_{\infty}(\lambda)$ (see Proposition \ref{prop6}), where 
$Z_k(\lambda) \in \mathcal{H}_{_{\mathcal{T},k}}$ is solution of~:
\begin{eqnarray*}
\left\{
\begin{array}{rcl}
(\nabla Z_k(\lambda),\nabla \xi) & = 0 & \forall \xi \in
H^1_{_{0,\mathcal{T}_k}}\\
\gamma Z_k(\lambda) & = \lambda & \mbox{ on } \Gamma~~~,
\end{array}
\right.
\end{eqnarray*}
and $Z_{\infty}(\lambda) \in \mathcal{H}_{_{\mathcal{T},\infty}}$ verifies~: 
\begin{eqnarray*}
\left\{
\begin{array}{ccc}
\Delta Z_{\infty}(\lambda) & = 0 & \mbox{ in } \Omega \\
\gamma Z_{\infty}(\lambda) & = \lambda & \mbox{ on } \Gamma~~~.
\end{array}
\right.
\end{eqnarray*}
\monitem Second, $Z_k(\lambda)$ is an approximation of the continuous harmonic function 
$Z_{\infty}(\lambda)$ on the mesh $\mathcal{T}_k.$ C\'ea's lemma \cite{Ci} insures that~:
$$\parallel Z_{\infty}(\lambda) - Z_k(\lambda) \parallel_{_{1,\Omega}}~\leq~C~
\parallel \Pi_{_{\mathcal{T}_k}} Z_{\infty}(\lambda) - Z_{\infty}(\lambda) \parallel_{_{1,\Omega}}~~~.$$
We fix $\varepsilon > 0$, independent of $\mathcal{T} \in \mathcal{U}_{\sigma}$. 
We have extended classical result of \cite{CR72} to problems with less regularity 
in Proposition \ref{prop9}, as $Z_{\infty}(\lambda)$ only belongs to $H^{2-\varepsilon}(\Omega)$ 
(see imbedding (\ref{reg})). Then, we choose $0 < \eta < 1 - \varepsilon$ and $\theta \geq \varepsilon$ 
is given by relation (\ref{theta}). So, we obtain~:
$$\parallel  \Delta Z_k(\lambda) \parallel_{_{-1,\Omega}}
~\leq~C~\parallel Z_{\infty}(\lambda) - Z_k(\lambda) \parallel_{_{1,\Omega}}
~\leq~C~h^{1-\theta}_k~\parallel Z_{\infty}(\lambda) \parallel_{_{2-\varepsilon,\Omega}}~~~,$$
where $\displaystyle{h_k \equiv h_{_{\mathcal{T}_k}} = \frac{h_{_{\mathcal{T}}}}{2^k}}.$
Finally, as $\mathcal{H}_{_{\mathcal{T},\infty}}$ is finite-dimensional, all the norms are equivalent. 
In particular, the $H^{2-\varepsilon}-$norm is equivalent to the $L^2-$norm, 
the constant involved depending on space dimension thus on $h_{_{\mathcal{T}}}$~:
$$\parallel  \Delta Z_k(\lambda) \parallel_{_{-1,\Omega}}
~\leq~C~\parallel Z_{\infty}(\lambda) - Z_k(\lambda) \parallel_{_{1,\Omega}}
~\leq~C~h^{1-\theta}_k~\parallel Z_{\infty}(\lambda) \parallel_{_{0,\Omega}}~~~.$$
\monitem Third, we bound below $\parallel Z_k(\lambda) \parallel_{_{0,\Omega}}$~:
\begin{eqnarray*}
\parallel Z_k(\lambda) \parallel_{_{0,\Omega}}
& \geq & \parallel Z_{\infty}(\lambda) \parallel_{_{0,\Omega}}~-~
         \parallel Z_{\infty}(\lambda) - Z_k(\lambda) \parallel_{_{0,\Omega}} \\
& \geq & \parallel Z_{\infty}(\lambda) \parallel_{_{0,\Omega}}~-~
         \parallel Z_{\infty}(\lambda) - Z_k(\lambda) \parallel_{_{1,\Omega}} \\
& \geq & \parallel Z_{\infty}(\lambda) \parallel_{_{0,\Omega}}~-~
         C~h^{1-\theta}_k~\parallel Z_{\infty}(\lambda) \parallel_{_{0,\Omega}}~~~.
\end{eqnarray*}
\monitem Thus, if $h_k$ is sufficiently small or $k$ sufficiently large in order to have~: $1-C h^{1-\theta}_k$ 
positive, we have~:
\begin{eqnarray*}
\rho(\lambda)~=~\frac{\parallel \Delta Z_k(\lambda) \parallel_{_{-1,\Omega}}}
                     {\parallel Z_k(\lambda) \parallel_{_{0,\Omega}}}
& \leq & \frac{C h^{1-\theta}_k \parallel Z_{\infty}(\lambda)\parallel_{_{0,\Omega}}}
              {(1-C h^{1-\theta}_k) \parallel Z_{\infty}(\lambda) \parallel_{_{0,\Omega}}} \\
& \leq & \frac{C h_{_{\mathcal{T}}}^{1-\theta} (2^{-k})^{1-\theta}}
              {(1-C h_{_{\mathcal{T}}}^{1-\theta} (2^{-k})^{1-\theta})}~\leq~\frac{1}{\delta}~~~,
\end{eqnarray*}
which is independent of $\lambda$ when $k \longrightarrow +\infty$. $\hfill \blacksquare$\\

By now, $\delta$ is taken equal to $1$.

\begin{defi}{Discrete harmonic formulation of the Stokes problem.} \\
For all mesh $\mathcal{T} \in \mathcal{U}_{\sigma}$ and all integer $k \geq 0$, we set~:
$$H^{1,k}_{_{\mathcal{T}}}~=~H^1_{_{0,\mathcal{T}}} \oplus \mathcal{H}_{_{\mathcal{T}_k}}~~~.$$
The discrete harmonic formulation of the Stokes problem consists in finding~:
\begin{eqnarray}
\psi_{_{\mathcal{T},k}} \in H^1_{_{0,\mathcal{T}}}~,~\omega_{_{\mathcal{T},k}} \in H^{1,k}_{_{\mathcal{T}}} 
& \mathrm{such} & \mathrm{that~:} \label{E15} \\
(\omega_{_{\mathcal{T},k}},\varphi) - (\nabla \psi_{_{\mathcal{T},k}},\nabla \varphi) 
& = & 0~~~~~~~~~~~~~\forall \varphi \in H^{1,k}_{_{\mathcal{T}}} \label{E16} \\
(\nabla \omega_{_{\mathcal{T},k}}, \nabla \xi) 
& = & (\mathbf{f},\rotV~\xi)~~~\forall \xi \in H^1_{_{0,\mathcal{T}}} \label{E17}
\end{eqnarray}
\end{defi}

\begin{rem} \label{rem4}
Notice that $\psi$ and $\omega^0$ are discretized on the initial mesh $\mathcal{T}$ and not on $\mathcal{T}_k.$ 
Notice also that $H^1_{_{0,\mathcal{T}}} \subset H^1_{_{0,\mathcal{T}_k}}$ because meshes are embedded.
\end{rem}

This formulation coincides with (\ref{E9-11}) when $k = + \infty$, and with
(\ref{E5})-(\ref{E8}) when $k = 0$. The next proposition shows that if $k$ is chosen sufficiently 
large, system (\ref{E15})-(\ref{E17}) has a solution which depends continuously on the datum.  

\begin{prop}{Existence and uniqueness of a solution to (\ref{E15})-(\ref{E17}).} \label{prop10} \\
If $\mathbf{f} \in (L^2(\Omega))^2$, if $\mathcal{T} \in \mathcal{U}_{\sigma}$ 
and if $k \geq K(1,\mathcal{T})$, with $K$ defined in Theorem \ref{th2}, 
the problem (\ref{E15})-(\ref{E17}) has a unique solution
$(\psi_{_{\mathcal{T},k}},\omega_{_{\mathcal{T},k}})$ in the space $H^1_{_{0,\mathcal{T}}} \times
H^{1,k}_{_{\mathcal{T}}}$ which depends continuously on the datum $\mathbf{f}$~:
$$\exists~C~>~0~~,~~\parallel \omega_{_{\mathcal{T},k}} \parallel_{_M}
~+~\parallel \nabla \psi_{_{\mathcal{T},k}} \parallel_{_{0,\Omega}}
~\leq~C~\parallel \mathbf{f} \parallel_{_{0,\Omega}}~~~.$$
\end{prop}

\proof \\
The proof is not different from the one of Proposition \ref{prop3}. We first split 
functions of $H^{1}_{_{\mathcal{T},k}}$ in a function $\varphi^{0,k}$ of $H^1_{_{0,\mathcal{T}}}$ plus 
a function $\varphi^{\Delta,k}$ in $\mathcal{H}_{_{\mathcal{T},k}}$~:
$$\varphi~=~\varphi^{0,k}~+~\varphi^{\Delta,k}~~~.$$
Using that each function $\varphi$ in $\mathcal{H}_{_{\mathcal{T},k}}$ verifies~:
$$(\nabla \varphi,\nabla \xi)~=~0~~,~~\forall \xi \in H^1_{_{0,\mathcal{T}}} \subset H^1_{_{0,\mathcal{T}_k}}~~~,$$
we rewrite the system (\ref{E16})-(\ref{E17}) and obtain~:
\begin{eqnarray}
(\omega^{0,k}_{_{\mathcal{T}}},\varphi^{0,k}) + (\omega^{\Delta,k}_{_{\mathcal{T}}},\varphi^{0,k})
- (\nabla \psi_{_{\mathcal{T},k}},\nabla \varphi^{0,k}) \! \! & = & \! \! 0~~~~~~~~~~~\forall \varphi^{0,k} \in
H^1_{_{0,\mathcal{T}}} \label{E18}\\
(\omega^{0,k}_{_{\mathcal{T}}},\varphi^{\Delta,k}) +
(\omega^{\Delta,k}_{_{\mathcal{T}}},\varphi^{\Delta,k})
\! \! & = & \! \! 0~~~~~~~~~~~\forall \varphi^{\Delta,k} \in \mathcal{H}_{_{\mathcal{T},k}} \label{E19}\\
(\nabla \omega^{0,k}_{_{\mathcal{T}}}, \nabla \xi) \! \! & = & \! \! (\mathbf{f},\rotV~\xi)~\forall \xi \in
H^1_{_{0,\mathcal{T}}}.
\label{E20}
\end{eqnarray}
Notice that functions in $\mathcal{H}_{_{\mathcal{T},k}}$ are not harmonic, so 
the $M$-norm is no more the $L^2-$norm. We need Theorem \ref{th2} to conclude that,  
if $k \geq K(1,\mathcal{T})$, we have~:
$$\forall \varphi \in \mathcal{H}_{_{\mathcal{T},k}}~~,~~\parallel \varphi \parallel_{_M}
~\geq~\sqrt{2}~\parallel \Delta \varphi \parallel_{_{-1,\Omega}}~~~,$$
which implies that the $L^2-$product $(\omega^{\Delta}_{_{\mathcal{T}}},\varphi^{\Delta})$ is elliptic on 
$\mathcal{H}_{_{\mathcal{T},k}}$ for the $M$-norm. Then, the same arguments as in the proof of Proposition 
\ref{prop3} can be used here, and we obtain the expected inequality. $\hfill \blacksquare$

\begin{rem}
Do not confuse~:
$$H^{1,k}_{_{\mathcal{T}}} \ni \varphi~=~\underbrace{\varphi^0}_{\in H^1_0(\Omega)}
                                      ~+~\underbrace{\varphi^{\Delta}}_{\in \mathcal{H}_{_{\mathcal{T},\infty}}(\Omega)}
                                      ~=~\underbrace{\varphi^{0,k}}_{\in H^1_{_{0,\mathcal{T}}}}
				      ~+~\underbrace{\varphi^{\Delta,k}}_{\in \mathcal{H}_{_{\mathcal{T},k}}}~~~.$$
\end{rem}

\subsection{Interpolation error}

\begin{prop}{Interpolation error.} \label{prop12}\\
Under the same assumptions as in Proposition \ref{prop5}, there exists an interpolation operator
$\phi^k_{_{\mathcal{T}}}$ from $\mathcal{H}(\Omega) \cap H^2(\Omega)$ to $\mathcal{H}_{_{\mathcal{T},k}}$ 
such that~:
$$\parallel \varphi - \phi^k_{_{\mathcal{T}}} \varphi \parallel_{_M}
~\leq~C~h^{1-\theta}_{_{\mathcal{T}}}~\parallel \varphi \parallel_{_{2,\Omega}}~~~,$$
where $\theta$ is an arbitrary real number such that $0 < \theta < 1$ and $C$ a constant independent 
of $k$ but dependent on $\theta$.
\end{prop}

\begin{defi} \label{rem7} {$~$}\\
Let define $\mathcal{P}^k_{_{\mathcal{T}}} : M(\Omega) \cap H^2(\Omega) = H^2(\Omega) \longrightarrow 
H^{1,k}_{_{\mathcal{T}}} = H^1_{_{0,\mathcal{T}}} \oplus \mathcal{H}_{_{\mathcal{T},k}}$~:
$$\mathcal{P}^k_{_{\mathcal{T}}} : H^2(\Omega) \ni \varphi = \varphi^0 + \varphi^{\Delta} 
\longmapsto \mathcal{P}^k_{_{\mathcal{T}}} \varphi = \Pi_{_{\mathcal{T}}} \varphi^0 
+ \phi^k_{_{\mathcal{T}}} \varphi^{\Delta} \in H^{1,k}_{_{\mathcal{T}}}~~~.$$
\end{defi}

\proof \\
We have supposed that $\varphi$ belongs to $H^2(\Omega)$. From Proposition \ref{prop5}, 
$\zeta = \phi^{\infty}_{_{\mathcal{T}}} \varphi^{\Delta}$ is an approximation of $\varphi^{\Delta}$ 
and verifies by definition~:
\begin{eqnarray*}
\left\{
\begin{array}{rcll}
\Delta \zeta & = & 0 & \mbox{ in } \Omega \\
\gamma \zeta & = & \Pi^{\Gamma}_{_{\mathcal{T}}}(\gamma \varphi) & \mbox{ on } \Gamma~~~.
\end{array}
\right.
\end{eqnarray*}
As $\lambda \in H^{3/2-\varepsilon}(\Gamma)$, for all $\varepsilon > 0$, and as $\Omega$ is a convex
polygon, $\zeta$ belongs to $H^{2-\varepsilon}(\Omega)$ \cite{ADN}, \cite{Gri}. 
We can interpolate $\zeta$ on the $k$ times refined mesh $\mathcal{T}_k$
because $H^{2-\varepsilon}(\Omega) \subset \mathcal{C}^0(\overline{\Omega})$. We choose 
$\varepsilon \leq \theta$ in order to satisfy (\ref{theta}) and we obtain~:
$$\parallel \zeta - \Pi_{_{\mathcal{T}_k}} \zeta \parallel_{_{1,\Omega}}
~\leq~C~h^{1-\theta}_k~\parallel \zeta \parallel_{_{2-\varepsilon,\Omega}}~~~.$$
So, we have successively~:
\begin{eqnarray*}
\parallel \zeta - \Pi_{_{\mathcal{T}_k}} \zeta \parallel_{_M} 
& \leq & \parallel \zeta - \Pi_{_{\mathcal{T}_k}} \zeta \parallel_{_{1,\Omega}} 
         \mbox{ by Proposition } \ref{rem2}, \\
& \leq & C~h^{1-\theta}_k~\parallel \zeta \parallel_{_{2-\varepsilon,\Omega}} \\
& \leq & C~h^{1-\theta}_k~\parallel \lambda \parallel_{_{3/2-\varepsilon,\Gamma}} 
         \mbox{ by continuity of the lifting,}\\
& \leq & C~h^{1-\theta}_k~\parallel \gamma \varphi \parallel_{_{3/2-\varepsilon,\Gamma}} \\
& \leq & C~h^{1-\theta}_k~\parallel \gamma \varphi \parallel_{_{3/2,\Gamma}} 
         \mbox{ by continuity of the trace operator } \\
&      & \hspace{4cm} \mbox{ from } H^{3/2-\varepsilon}(\Gamma) \mbox{ to }  H^{3/2}(\Gamma),\\
& \leq & C~h^{1-\theta}_k~\parallel \varphi \parallel_{_{2,\Omega}}
         \mbox{ by continuity of the trace.}
\end{eqnarray*}
Thus, we define $\phi^k_{_{\mathcal{T}}} \varphi~=~\Pi_{_{\mathcal{T}_k}} \zeta~=~
\Pi_{_{\mathcal{T}_k}}(\phi^{\infty}_{_{\mathcal{T}}} \varphi)$ and we obtain~:
\begin{eqnarray*}
\parallel \varphi - \phi^k_{_{\mathcal{T}}} \varphi \parallel_{_M} 
& =    & \parallel \varphi - \Pi_{_{\mathcal{T}_k}} \zeta \parallel_{_M} \\
& \leq & \parallel \varphi - \phi^{\infty}_{_{\mathcal{T}}} \varphi \parallel_{_M}~+~
         \parallel \underbrace{\phi^{\infty}_{_{\mathcal{T}}} \varphi}_{= \zeta} 
	 - \underbrace{\Pi_{_{\mathcal{T}_k}}(\phi^{\infty}_{_{\mathcal{T}}} \varphi)}
	             _{= \Pi_{_{\mathcal{T}_k}} \zeta} \parallel_{_M} \\
& \leq & C~h_{_{\mathcal{T}}}~\parallel \varphi \parallel_{_{2,\Omega}}
         ~+~\parallel \zeta - \Pi_{_{\mathcal{T}_k}} \zeta \parallel_{_M}
         \mbox{ by Proposition \ref{prop5}}\\
& \leq & C~h_{_{\mathcal{T}}}~\parallel \varphi \parallel_{_{2,\Omega}}
         ~+~C~h^{1-\theta}_k~\parallel \varphi \parallel_{_{2,\Omega}} \\
& \leq & C~(h_{_{\mathcal{T}}}~+~h^{1-\theta}_{_{\mathcal{T}}} (2^{-k})^{1-\theta}) 
         \parallel \varphi \parallel_{_{2,\Omega}}~~~.
\end{eqnarray*}
As $h_k~\leq~h_{_{\mathcal{T}}}$, Proposition \ref{prop12} is proved. $\hfill \blacksquare$

\begin{rem} \label{rem8}
Notice that the interpolation error is independent of $k$ when $k$ is large enough,  
because $h_k \leq h_{_{\mathcal{T}}}$. If we need $k$ sufficiently large, 
it is for existence (Proposition \ref{prop10}) and stability (Proposition \ref{prop11}). 
\end{rem}

\subsection{Error estimates}

\begin{prop}{Stability of problem (\ref{E15})-(\ref{E17}).} \label{prop11} \\
Let $\delta = 1$, $\mathcal{T}$ be a mesh of the regular family $\mathcal{U}_{\sigma}$, and
$K(1,\mathcal{T})$ chosen as in Theorem \ref{th2}. Let 
$(\eta_{_{\mathcal{T},k}},\theta_{_{\mathcal{T},k}}) \in H^1_{_{0,\mathcal{T}}} 
\times H^{1,k}_{_{\mathcal{T}}}$ be the solution of the auxiliary problem~:
\begin{eqnarray} \label{P2st}
\left\{
\begin{array}{rcll}
(\theta_{_{\mathcal{T},k}},\varphi) - (\nabla \eta_{_{\mathcal{T},k}},\nabla \varphi)
& = & (g,\varphi) + \langle \Delta \varphi,m\rangle_{_{-1,1}} 
& \forall \varphi \in H^{1,k}_{_{\mathcal{T}}}\\
(\nabla \theta_{_{\mathcal{T},k}},\nabla \xi)
& = & \langle l,\xi\rangle_{_{-1,1}} 
& \forall \xi \in H^1_{_{0,\mathcal{T}}}
\end{array}
\right.
\end{eqnarray}
where $g = \mathcal{P}^k_{_{\mathcal{T}}} \omega - \omega \in L^2(\Omega)$, $m = \Pi_{_{\mathcal{T}}} \psi -
\psi \in H^1_0(\Omega)$ and $l = \Delta (\omega - \mathcal{P}^k_{_{\mathcal{T}}} \omega) \in H^{-1}(\Omega)$.  
$\mathcal{P}^k_{_{\mathcal{T}}}$ is the interpolation operator from $M(\Omega)$ to $H^{1,k}_{_{\mathcal{T}}}$, 
defined in Definition \ref{rem7} and $\Pi_{_{\mathcal{T}}}$ the classical interpolation operator on $\mathcal{T}$. 
Then, if $k \geq K(1,\mathcal{T})$, there exists a constant $C$ independent of $\sigma$ and $k$ such that~: 
$$\parallel \theta_{_{\mathcal{T},k}} \parallel_{_M}~+~\parallel \nabla \eta_{_{\mathcal{T},k}} \parallel_{_{0,\Omega}}
~\leq~C~\left( \parallel g \parallel_{_{0,\Omega}} + \parallel \nabla m \parallel_{_{0,\Omega}} + 
\parallel l \parallel_{_{-1,\Omega}} \right)~~~.$$
\end{prop}

\proof \\
The auxiliary problem is obtained by the same techniques used in Proposition \ref{ieme} and 
arguments are identical to those in Proposition \ref{prop4}.
But, for $\theta^{\Delta,k}_{_{\mathcal{T}}} \in \mathcal{H}_{_{\mathcal{T},k}}$, 
the $L^2$-norm is no more the $M$-norm. We need Theorem \ref{th2} to conclude that
if $k \geq K(1,\mathcal{T})$, we have~:
$$\parallel \theta^{\Delta,k}_{_{\mathcal{T}}} \parallel_{_{0,\Omega}}
~\geq~\parallel \Delta \theta^{\Delta,k}_{_{\mathcal{T}}} \parallel_{_{-1,\Omega}}~~~.$$
$\hfill \blacksquare$\\

The previous stability (Proposition \ref{prop11}) and interpolation errors (relation (\ref{phi0}) 
for the $H^1_0-$part of the vorticity and for the stream function; Proposition \ref{prop12} for 
the harmonic part) gives the following result~:
 
\begin{theorem}{Convergence result.}\label{th3}\\
Problem (\ref{E16}-\ref{E17}) has a unique solution
$(\psi_{_{\mathcal{T},k}},\omega_{_{\mathcal{T},k}}) \in H^1_{_{0,\mathcal{T}}}\times H^{1,k}_{_{\mathcal{T}}}$.
If $(\psi,\omega)$, belonging to $(H^2(\Omega) \cap H^1_0(\Omega)) \times H^2(\Omega)$, are solution of 
the continuous problem (\ref{P}), we have~:
$$\forall~\theta~\in~]0,1[~,~\exists~C~=~C(\theta)~>~0~,~\forall~\mathcal{T}~\in~\mathcal{U}_{\sigma}~,
~\exists~K(1,\mathcal{T})~,~\mbox{if}~k~\geq~K(1,\mathcal{T})$$
$$\parallel \omega - \omega_{_{\mathcal{T},k}}  \parallel_{_M}
   ~+~\parallel \psi - \psi_{_{\mathcal{T},k}} \parallel_{_{1,\Omega}}~\leq~C~h^{1-\theta}_{_{\mathcal{T}}}~
   (\parallel \omega \parallel_{_{2,\Omega}} + \parallel \psi \parallel_{_{2,\Omega}})~~~.$$
\end{theorem}

\section{Numerical results}

Numerical experiments have been performed on a unit square with an analytical solution (test of 
Bercovier-Engelman \cite{BE}). 
\begin{eqnarray*}
f_1(x,y) & = & 256 \left( x^2 (x-1)^2 (12y-6) + y (y-1) (2y-1) (12x^2 -12x+2) \right) \\
& + & (y-0.5) \\
f_2(x,y) & = & - f_1(y,x) \\
\psi(x,y) & = & -128 (y^2 (y-1)^2 x^2 (x-1)^2) \\
\omega(x,y) & = & 256 (y^2 (y-1)^2 (6x^2-6x+1) + x^2 (x-1)^2 (6y^2-6y+1)).   
\end{eqnarray*}
As it is said in Theorem \ref{th2}, $K$ depends on the mesh, so we have worked with three 
different unstructured meshes obtained with EMC2, mesh generator of Modulef (Bernadou and al. \cite{Ber}). 
Notice that structured meshes give good results without any 
stabilization (see \cite{DSS98}, \cite{GR}). The uncoupled algorithm is the following~:
$$
\begin{array}{lrcll}
\mbox{Find } \omega^{0,k}_{_{\mathcal{T}}} \in H^1_{_{0,\mathcal{T}}} \mbox{ such that} &
(\nabla \omega^{0,k}_{_{\mathcal{T}}},\nabla \xi) & = & (\mathbf{f},\mathbf{rot}~\xi)
& \forall \xi \in H^1_{_{0,\mathcal{T}}} \\
\mbox{Find } \omega^{\Delta,k}_{_{\mathcal{T},k}} \in \mathcal{H}_{_{\mathcal{T}_k}} \mbox{ such that} &
(\omega^{\Delta,k}_{_{\mathcal{T},k}},\chi)       & = & - (\omega^{0,k}_{_{\mathcal{T}}},\chi)
& \forall \chi \in \mathcal{H}_{_{\mathcal{T}_k}} \\
\mbox{Find } \psi_{_{\mathcal{T},k}} \in H^1_{_{0,\mathcal{T}}} \mbox{ such that} &
(\nabla \psi_{_{\mathcal{T},k}},\nabla \xi)       & = & (\omega^{0,k}_{_{\mathcal{T}}} + 
\omega^{\Delta,k}_{_{\mathcal{T},k}},\xi)
& \forall \xi \in H^1_{_{0,\mathcal{T}}}~.
\end{array}
$$
An advantage of this algorithm is to avoid the well-known difficulties related to the 
resolution of mixed linear systems (see \textit{e.g.} \cite{BDM85}, \cite{ES96}).

\begin{rem}
It is clear that $\omega^{0,k}_{_{\mathcal{T}}}$ does not depend on $k$.
\end{rem}

The longer and expensive step is the second one. By the way, we have to calculate {\bf{all}} the 
discrete harmonic functions (as many as edges of the initial triangulation on the boundary) 
{\bf{on the refined meshes}}. We present the three different initial meshes used in those examples 
Figures \ref{f:meshab}, \ref{f:meshc} and the number of vertices of each refined meshes (Table \ref{tab:0}).

\begin{figure}    [H]  \centering
\centerline  {\includegraphics[width=.74\textwidth]   {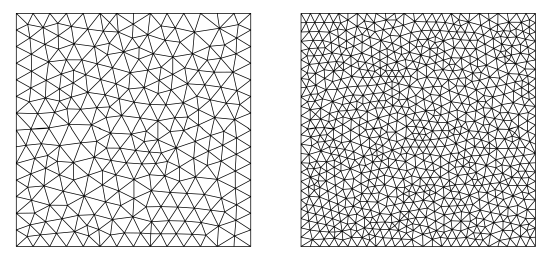}}
\caption{Meshes A and B}
\label{f:meshab}
\end{figure}
%
\begin{figure}    [H]  \centering
\centerline  {\includegraphics[width=.37\textwidth]   {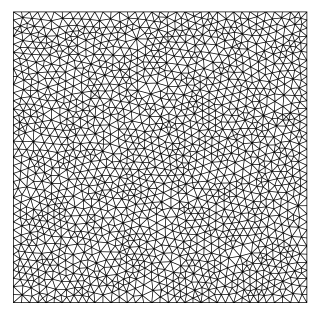}}
\caption{Mesh C}
\label{f:meshc}
\end{figure}

\begin{table}[hbt]
\begin{center}
\begin{tabular} {|l|c|c|c|c|c|}
\hline 
{Mesh } & {$k=0$} & {$k=1$} & {$k=2$} & {$k=3$} & {$k=4$} \\
\hline  
A & $270$ & $1021$ & $3969$ & $15649$ & $62145$ \\ \hline 
B & $832$ & $3225$ & $12697$ & $50385$ & $200737$ \\ \hline 
C & $1667$ & $6521$ & $25793$ & $102593$ & $409217$ \\ \hline 
\end{tabular}
\end{center}
\caption{Number of vertices of the k times refined meshes.}
\label{tab:0}
\end{table}

The analytical vorticity attains its extremum on the middle of each edge of the square and 
its value is then $+16.00$ (see Figure \ref{f:tourbE}). In fact, without any stabilization, 
extrema of the vorticity exploded on the boundary. Table \ref{tab:1} gives extrema of the 
vorticity obtained after $0$ (corresponds to the non-stabilized scheme), $1, 2, 3$ and $4$ 
refinements of the mesh. The Figure \ref{f:tourbE} gives the values of the vorticity along the boundary. 
Let us recall that the number of discrete harmonic functions is equal to the number of vertices on the boundary.

\begin{figure}    [H]  \centering
\centerline  {\includegraphics[width=.80\textwidth]   {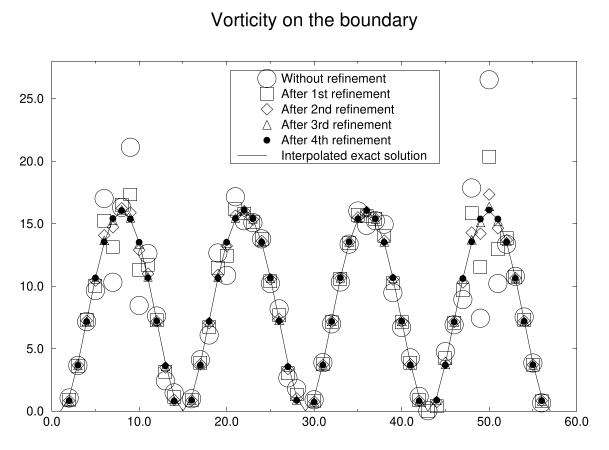}}
\caption{Value of the vorticity on the boundary (mesh A).}
\label{f:tourbE}
\end{figure}


\begin{table}[hbt]
\begin{center}
\begin{tabular} {|l|c|c|c|c|c|c|}
\hline 
{Mesh } & {Number of} & {$k=0$} & {$1$} & {$2$} & {$3$} & {$4$} \\ 
{} & {discrete} & {} & {} & {} & {} & {} \\ 
{} & {harmonics} & {} & {} & {} & {} & {} \\ \hline
A $: 482$  elements &  56 & $26.52$ & $20.37$ & $17.34$ & $16.37$ & $16.10$ \\ \hline 
B $: 1562$ elements & 100 & $27.08$ & $20.78$ & $17.41$ & $16.34$ & $16.05$ \\ \hline 
C $: 3188$ elements & 144 & $31.77$ & $22.94$ & $18.20$ & $16.60$ & $16.15$ \\ \hline 
\end{tabular}
\end{center}
\caption{Extremum of the vorticity after $k$ refinements.}
\label{tab:1}
\end{table}


We observe that $4$ refinements seem sufficient to obtain satisfying results. 
This number of necessary refinements is confirmed by Figures \ref{meshE}, 
\ref{meshF}, \ref{meshG}, which give the relative error on the two fields (stream function 
and vorticity) versus the number of refinements. Error is no more reduced between fourth 
and fifth refinement, so we can estimate that $k = 4$ shall be enough and that stability 
of the scheme is reached.\\ 

\begin{figure}    [H]  \centering
\centerline  {\includegraphics[width=.85\textwidth]   {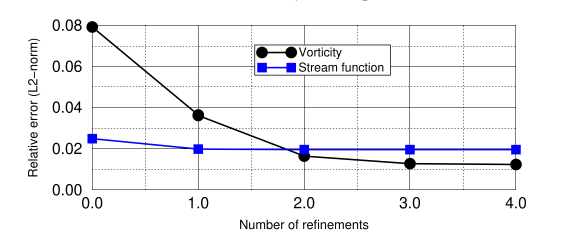}}
\caption{Relative error in $L^2-$norm versus mesh size on mesh A.}
\label{meshE}
\end{figure}
%
\begin{figure}    [H]  \centering
\centerline  {\includegraphics[width=.85\textwidth]   {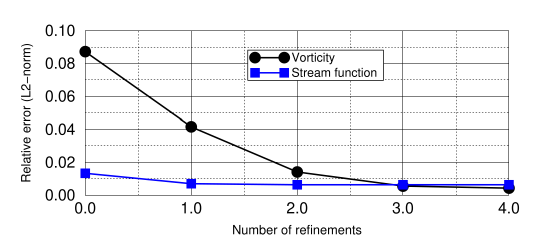}}
\caption{Relative error in $L^2-$norm versus mesh size on mesh B.}
\label{meshF}
\end{figure}
%
\begin{figure}    [H]  \centering
\centerline  {\includegraphics[width=.85\textwidth]   {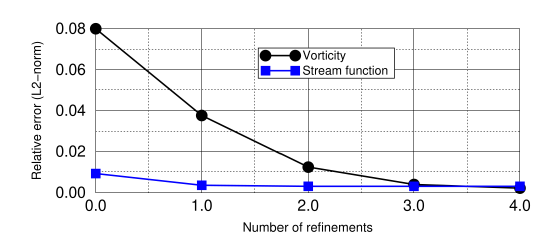}}
\caption{Relative error in $L^2-$norm versus mesh size on mesh C.}
\label{meshG}
\end{figure}
%
\begin{figure}    [H]  \centering
\centerline  {\includegraphics[width=.85\textwidth]   {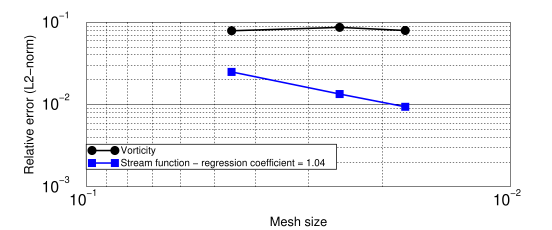}}
\caption{Convergence order without stabilization.}
\label{convns}
\end{figure}
%
\begin{figure}    [H]  \centering
\centerline  {\includegraphics[width=.85\textwidth]   {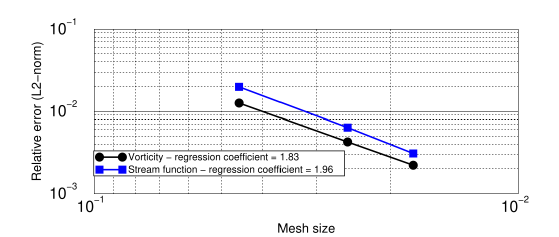}}
\caption{Convergence order after $4$ refinements.}
\label{convst}
\end{figure}

We recall that problem (\ref{E5})-(\ref{E8}) is not stable and error bounds are in 
$h^{1/2}_{_{\mathcal{T}}}$ for the $L^2-$norm of the vorticity and $h^{1-\varepsilon}_{_{\mathcal{T}}}$ 
for the $H^1-$norm of the stream function (see Figure \ref{convns} or \cite{GR}).
With techniques described in this paper, we obtain a convergence of order 1 (see Theorem \ref{th3}).
Notice that numerical results, given on Figure \ref{convst}, are better than expected~: 
almost $\mathcal{O}(h^2_{_{\mathcal{T}}})$ for the $L^2-$norm of the vorticity and of the stream function !

\section{Conclusion}

We have studied the well-posed Stokes problem in stream function and vorticity
formulation. We have shown that using a space of {\bf{real}} harmonic functions is sufficient to obtain 
convergence with better estimations than the previous ones. Then, we have proposed a way of approaching 
numerically the space of real harmonic functions by refining the initial mesh. We have shown 
theoretically that, by this way, we keep convergence with order $1$. Numerically, it seems that 
convergence is obtained with order $2$ on the $L^2-$norm of the vorticity. \\

We are actually studying how to prove the convergence with order $2$ and how to use the same technique for 
the vorticity-velocity-pressure formulation of the Stokes problem, introduced in \cite{Du92}, which is an 
alternative in three-dimensional domain to recent work of Amara and al. \cite{ABD}. \\

\noindent {\bf{Acknowledgements}}\\
The authors thank Beno\^{\i}t Desjardins and Emmanuel Grenier for helpful discussions during summer 1998, and 
the referees for their precise remarks and suggestions that have been introduced in this corrected version.





\end{document}